\documentclass[12pt]{article}

\usepackage[margin=1in]{geometry}  
\usepackage{graphicx}              
\usepackage{amsmath}               
\usepackage{amsfonts}              
\usepackage{amsthm}                
\usepackage{amssymb}

\begin{document}

\nocite{*}

\title{Tetrahedral quartics in Projective Space.}

\author{Evgeny Mayanskiy}

\maketitle

\begin{abstract}
  We study tetrahedral quartics in projective space. We address their projective geometry, Neron-Severi lattice and automorphism group. 
\end{abstract}

\setcounter{section}{-1}
\section{Introduction}

\subsection{Definition}

{\bf Definition.} A {\it tetrahedral surface in } ${\mathbb P}^3$ is a surface in ${\mathbb P}^3$, containing edges of a tetrahedron.\\ 

\begin{center}
  \includegraphics{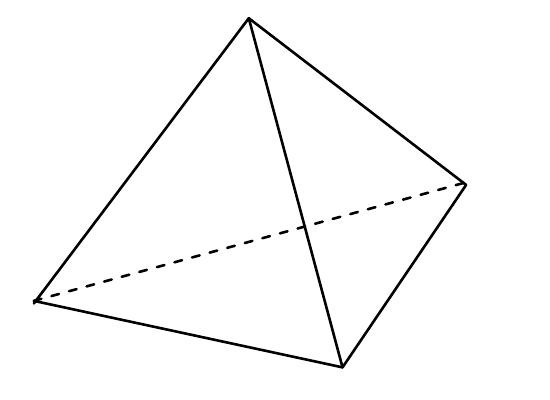}
\end{center}

{\bf Example.} There is a unique irreducible tetrahedral cubic surface in ${\mathbb P}^3$ - the cubic surface with $4$ nodes.\\

We will focus on the next simplest case - tetrahedral quartics in ${\mathbb P}^3$.

{\bf Lemma 0.1.} {\it A tetrahedral surface in ${\mathbb P}^3$ of degree $d\geq 4$ has equation $F(X_0,X_1,X_2,X_3)=0$ with polynomial 
$$
F(X_0,X_1,X_2,X_3)=A(X_0,X_1,X_2)\cdot X_0X_1X_2+B(X_0,X_1,X_3)\cdot X_0X_1X_3+
$$
$$
C(X_0,X_2,X_3)\cdot X_0X_2X_3+D(X_1,X_2,X_3)\cdot X_1X_2X_3+{\delta}\cdot X_0X_1X_2X_3,
$$
where $A, B, C, D$ are homogeneous polynomials of degree $d-3$ and $\delta$ is a homogeneous polynomial of degree $d-4$.}\\

{\it Proof:} Elementary. {\it QED}\\

{\bf Corollary 0.2.} {\it A tetrahedral quartic $X$ in ${\mathbb P}^3$ has equation $F(X_0,X_1,X_2,X_3)=0$ with polynomial 
$$
F(X_0,X_1,X_2,X_3)=A(X_0,X_1,X_2)\cdot X_0X_1X_2+B(X_0,X_1,X_3)\cdot X_0X_1X_3+
$$
$$
C(X_0,X_2,X_3)\cdot X_0X_2X_3+D(X_1,X_2,X_3)\cdot X_1X_2X_3+{\delta}\cdot X_0X_1X_2X_3,
$$
where $A(X_0,X_1,X_2)=a_0X_0+a_1X_1+a_2X_2, B(X_0,X_1,X_3)=b_0X_0+b_1X_1+b_3X_3, C(X_0,X_2,X_3)=c_0X_0+c_2X_2+c_3X_3, D(X_1,X_2,X_3)=d_1X_1+d_2X_2+d_3X_3$ are linear forms and $\delta \in \mathbb C$ is a constant. The linear forms $A, B, C, D$ represent residual lines of intersection of $X$ with the faces of the tetrahedron.}\\

\begin{center}
  \includegraphics{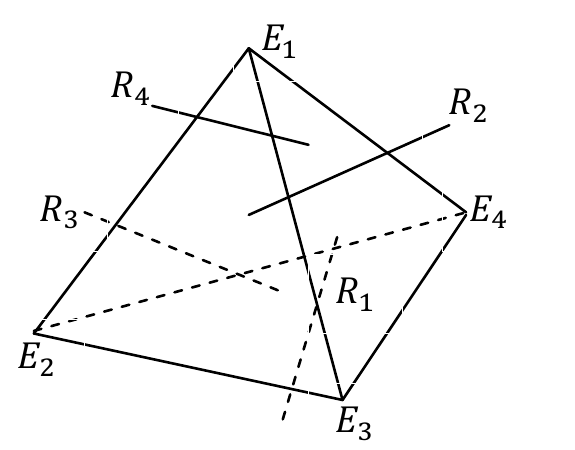}
\end{center}

In what follows we restrict our attention to tetrahedral quartics with $\delta = 1$.\\

{\bf Remark.} Tetrahedral quartics are parametrized by ${\mathbb C}^{12}$ (the coefficients of the linear forms $A, B, C, D$) or more precisely by the quotient ${{\mathbb C}^{12}}/ ( {{\mathbb C}^{*}}^{4}/{\mathbb C}^{*} )$, which exists as a quasiprojective variety and can be constructed, for example, using GIT \cite{Mumford}. One can describe easily stable and semi-stable tetrahedral quartics (using, for example, a general description of stability for torus actions in \cite{DolgachevGIT}).\\  

In what follows we will call a tetrahedral quartic 'general', if it corresponds to a point in a Zariski open subset $U \subset {\mathbb C}^{12}$, and 'very general', if it corresponds to a point in a subset $U^{an}\subset {\mathbb C}^{12}$ open in the analytic topology, which is a complement of a countable union of algebraic hypersurfaces in ${\mathbb C}^{12}$.\\

{\bf Corollary 0.3.} {\it A general tetrahedral quartic has $4$ nodes (vertices of the tetrahedron) and no other singularities.}\\

{\bf Corollary 0.4.} {\it For a general tetrahedral quartic $X$ projection $\phi \colon X\rightarrow {\mathbb P}^2$ from any of the $4$ nodes onto the opposite face of the tetrahedron represents $X$ as a (blow-up/ blow-down of a ) double cover of ${\mathbb P}^2$ ramified over an irreducible sextic $C(X)$ in ${\mathbb P}^2$ having:
\begin{itemize}
\item $3$ cusps (which are vertices of a triangle, each of whose sides touches $C(X)$ at a smooth point) and no other singularities,
\item one tritangent line $L$,
\item one tritangent conic $Q$ passing through the cusps (i.e. a smooth conic $Q$ in ${\mathbb P}^2$, which passes through the $3$ cusps of $C(X)$ and is tangent to $C(X)$ at $3$ smooth points),
\item a smooth cubic passing through the $3$ cusps of $C(X)$, through the $3$ points of tangency of $C(X)$ with the sides of the triangle at smooth points of $C(X)$, through the $3$ points of tangency of $C(X)$ with $L$ and through the $3$ points of tangency of $C(X)$ with $Q$.
\end{itemize}}

\begin{center}
  \includegraphics{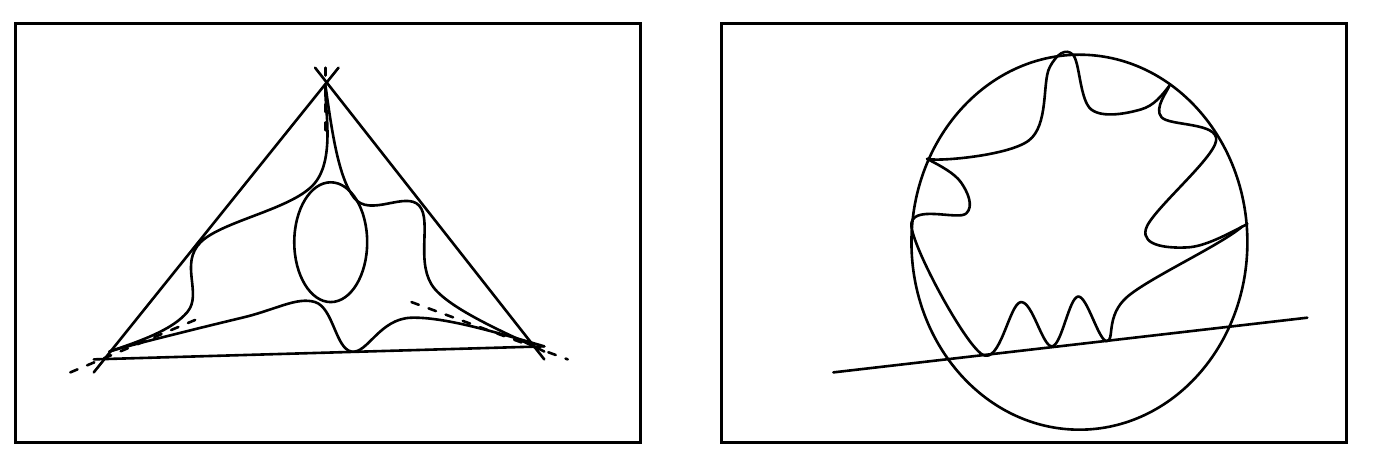}
\end{center}

{\bf Remark.} We will see later that the plane curve $C(X)$ characterizes tetrahedral quartic $X$. It is also clear that any irreducible sextic in ${\mathbb P}^2$ with $3$ cusps (which form a triangle whose sides touch the sextic at $3$ other points), tritangent line (which does not pass through the cusps and the points of tangency of the sextic with the sides of the triangle), tritangent smooth conic (which passes through the cusps, which is not tangent to the sextic at its cusps and which does not pass through its points of tangency with the triangle and the line) and a smooth cubic (passing through the $12$ points of tangency described above) appears as the branch curve of the projection of a tetrahedral quartic in ${\mathbb P}^3$ from one of its nodes.\\ 

\subsection{Special cases}

Tetrahedral quartics admit a plethora of interesting degenerations. See, for example, \cite{Bauer} (or Art. 116 of \cite{Jessop}). Let us notice here that as a special case of tetrahedral quartics (which appears when the residual lines $R_1, R_2, R_3, R_4$ are coplanar) one obtains hessian quartics. Another nice special case occurs, when the residual lines lie on a quadric. These can be naturally gathered in groups of $16$ in general (corresponding to the fact that a smooth quadric in ${\mathbb P}^3$ intersects each face of a tetrahedron in $2$ lines in this case). We will consider these 'quadratic perturbations' of hessians elsewhere.\\  

\subsection{Notation}

We will use the following notation for some natural curves on a general tetrahedral quartic $X$ and its minimal resolution of singularities $\tilde{X}$.\\

$E_i$, where $i=1,2,3,4$ will be a node of $X$ (or the corresponding exceptional divisor on $\tilde{X}$).\\

\begin{center}
  \includegraphics{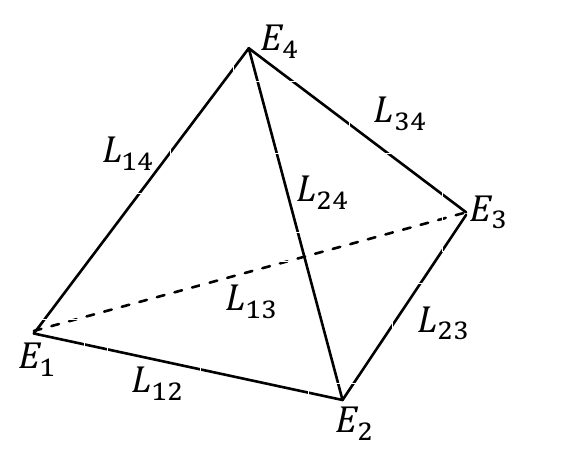}
\end{center}

$L_{ij}$, where $1\leq i < j \leq 4$ will be a line in ${\mathbb P}^3$ (lying on $X$) connecting nodes $E_i$ and $E_j$ (or its strict transform on $\tilde{X}$).\\

\begin{center}
  \includegraphics{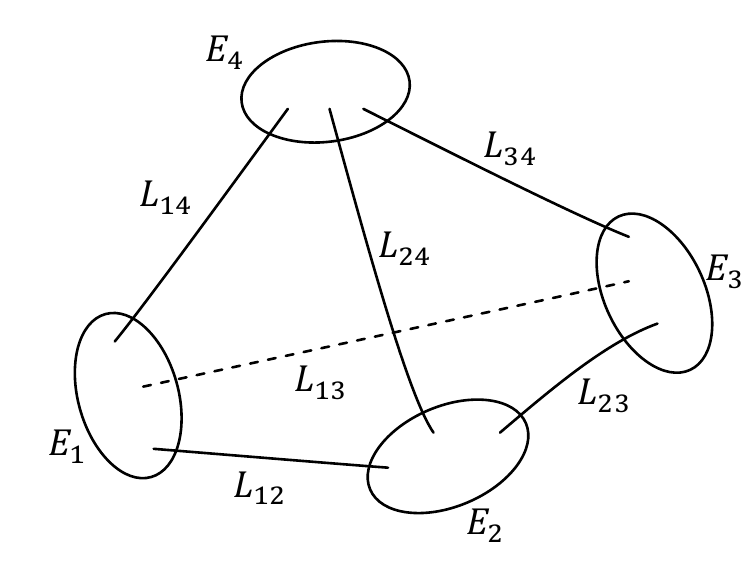}
\end{center}

$R_{j}$, where $j=1,2,3,4$ will be a residual line on $X$ (or its strict transform on $\tilde{X}$) which is opposite to the node $E_j$ .\\

\begin{center}
  \includegraphics{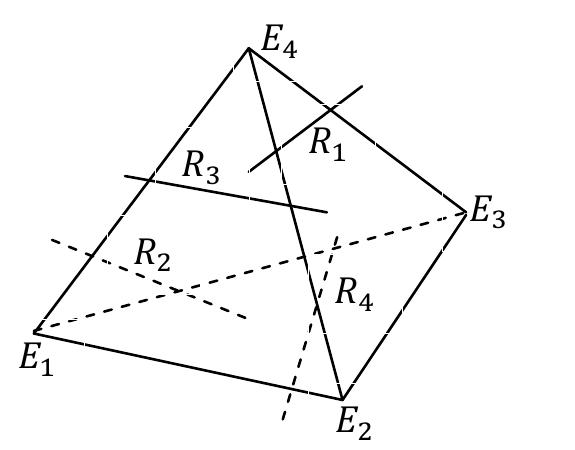}
\end{center}
We will also use the following notation for divisor classes on $\tilde{X}$: $H=L_{12}+L_{24}+L_{14}+E_1+E_2+E_4+R_3=$ $L_{24}+L_{23}+L_{34}+E_2+E_3+E_4+R_1=$ $L_{13}+L_{14}+L_{34}+E_1+E_3+E_4+R_2=$ $L_{12}+L_{13}+L_{23}+E_1+E_2+E_3+R_4\in NS(\tilde{X})$ is the class of the hyperplane section of $X$ in ${\mathbb P}^3$, $A=L_{12}+L_{13}+L_{14}+L_{23}+L_{24}+L_{34}+E_1+E_2+E_3+E_4+R_1+R_2+R_3+R_4\in NS(\tilde{X})$ and $A_0=3H-E_1-E_2-E_3-E_4\in NS(\tilde{X})$ will be shown to be ample divisor classes on $\tilde{X}$.\\

\subsection{Mirror Symmetry}

Let $\pi \colon \tilde{X}\rightarrow X$ be the minimal resolution of singularities. If $X$ is a general tetrahedral quartic in ${\mathbb P}^3$, then $\pi$ is just the blow-up of the $4$ nodes on $X$, and $\tilde{X}$ is a $K3$ surface.\\

General tetrahedral quartics come in 'mirror' pairs: $X_1 \subset {\mathbb P}^3$ and $X_2 \subset {\mathbb P}^3$ have isomorphic minimal resolutions of singularities $\tilde{X_1}\cong \tilde{X_2}$, but $X_1$ and $X_2$ are not projectively isomorphic.\\ 

\begin{center}
  \includegraphics{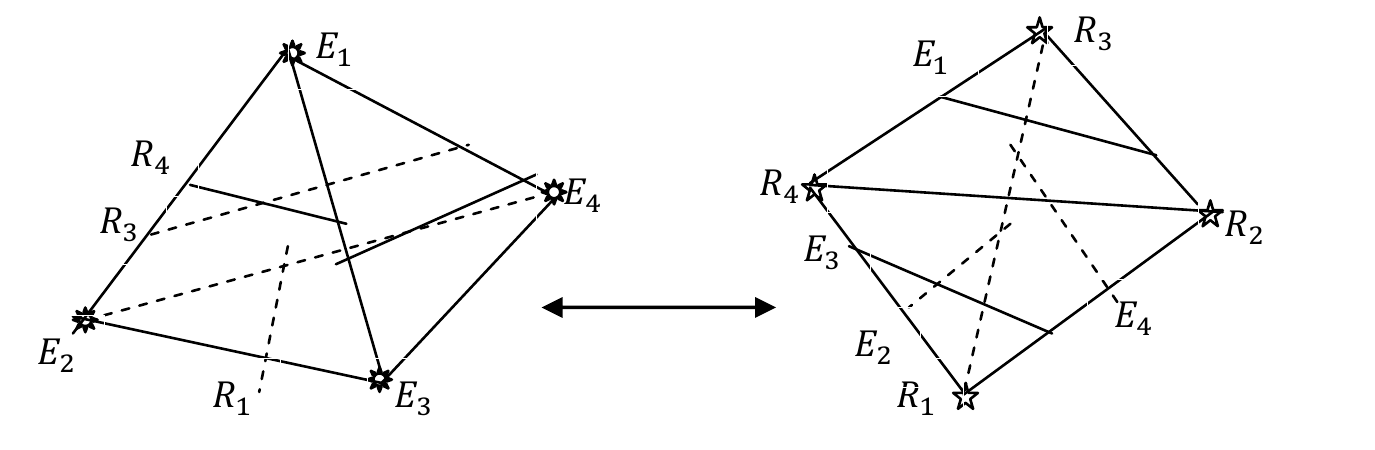}
\end{center}
Indeed, given a general tetrahedral quartic $X\subset {\mathbb P}^3$ one can interchange nodes and residual lines and obtain its mirror dual tetrahedral quartic $X^{v}\subset {\mathbb P}^3$. Explicitly, one can use the complete linear system $\mid 3\cdot H-2\cdot (E_1+E_2+E_3+E_4)-(L_{12}+L_{13}+L_{14}+L_{23}+L_{24}+L_{34}) \mid$ on $\tilde{X}$ , which is cut out on $X$ by the cubic surfaces in ${\mathbb P}^3$ passing through the edges of the tetrahedron. It is immediate that $(X^{v})^{v}=X$.\\ 

The fact that $X$ and $X^{v}$ are not projectively isomorphic follows immediately from the comparison of cross-ratios of the points of intersection of the residual lines and nodes on a general tetrahedral quartic with the edges of the tetrahedron.\\

Let $X\subset {\mathbb P}^3$ be a general tetrahedral quartic with equation $F(X_0,X_1,X_2,X_3)=0$ with polynomial $F(X_0,X_1,X_2,X_3)=A(X_0,X_1,X_2)\cdot X_0X_1X_2+B(X_0,X_1,X_3)\cdot X_0X_1X_3+C(X_0,X_2,X_3)\cdot X_0X_2X_3+D(X_1,X_2,X_3)\cdot X_1X_2X_3+ X_0X_1X_2X_3$, where $A(X_0,X_1,X_2)=a_0X_0+a_1X_1+a_2X_2, B(X_0,X_1,X_3)=b_0X_0+b_1X_1+b_3X_3, C(X_0,X_2,X_3)=c_0X_0+c_2X_2+c_3X_3, D(X_1,X_2,X_3)=d_1X_1+d_2X_2+d_3X_3$. Let ${\lambda}_{ij}$ be the cross-ratio on $L_{ij}$ of the nodes $E_i$, $E_j$ and of the points of intersection of $L_{ij}$ with the two residual lines. Choose notation in such a way that $E_1=(1:0:0:0), E_2=(0:1:0:0), E_3=(0:0:1:0), E_4=(0:0:0:1)$ in ${\mathbb P}^3$.\\

{\bf Lemma 0.5.} {\it Assume that $X\subset {\mathbb P}^3$ is general. Then ${\lambda}_{12}=\frac{a_1b_0}{a_0b_1}$, ${\lambda}_{13}=\frac{a_2c_0}{a_0c_2}$, ${\lambda}_{14}=\frac{b_3c_0}{b_0c_3}$, ${\lambda}_{23}=\frac{a_2d_1}{a_1d_2}$, ${\lambda}_{24}=\frac{b_3d_1}{b_1d_3}$, ${\lambda}_{34}=\frac{c_3d_2}{c_2d_3}$.\\
In particular, ${\lambda}_{12}, {\lambda}_{13}, {\lambda}_{14}, {\lambda}_{23}, {\lambda}_{24}, {\lambda}_{34}$ are general and independent except that they satisfy one relation: ${\lambda}_{12}\cdot \frac{1}{{\lambda}_{13}}\cdot {\lambda}_{14}\cdot {\lambda}_{23}\cdot \frac{1}{{\lambda}_{24}}\cdot {\lambda}_{34}=1$.}\\

{\bf Corollary 0.6.} {\it If $X\subset {\mathbb P}^3$ be a general tetrahedral quartic, then $X$ and $X^{v}$ are not projectively isomorphic.}\\

{\it Proof:} If $X\subset {\mathbb P}^3$ is general, then $X$ contains only one tetrahedron, which should transform to itself by a projective linear transformation. Hence a projective isomorphism which identifies $X$ and its mirror dual $X^{v}$ would transform a triangle $L_{12},L_{13},L_{23}$ in ${\mathbb P}^3$ into a triple of concurrent lines in ${\mathbb P}^3$, which is impossible. {\it QED}\\

{\bf Corollary 0.7.} {\it If an automorphism of the minimal resolution of singularities $\tilde{X}$ of a general tetrahedral quartic $X\subset {\mathbb P}^3$ permutes the curves $\{  L_{12}, L_{13}, L_{14}, L_{23}, L_{24}, L_{34}, $ $E_1, E_2, E_3, E_4, $ $R_1, R_2, R_3, R_4 \}$ among themselves, then it is the identity.}\\

{\it Proof:} Since cross-ratios ${\lambda}_{ij}$ are pairwise distinct, the permutation of the $14$ curves listed in the Corollary is trivial. Since each of the curves intersects the others at at least $3$ points, each of these curves is fixed pointwise. Hence such an automorphism of $\tilde{X}$ comes from an automorphism of ${\mathbb P}^3$, which fixes the tetrahedron pointwise. Since the residual lines are also pointwise fixed, such an automorphism should be identity. {\it QED}\\

\subsection{Elliptic fibrations}

Tetrahedral quartics admit several natural elliptic fibrations.\\

In particular, one has $6$ elliptic fibrations coming from pencils of planes in ${\mathbb P}^3$ containing the line $L_{ij}$, $1\leq i < j\leq 4$, and\\ 

\begin{center}
  \includegraphics{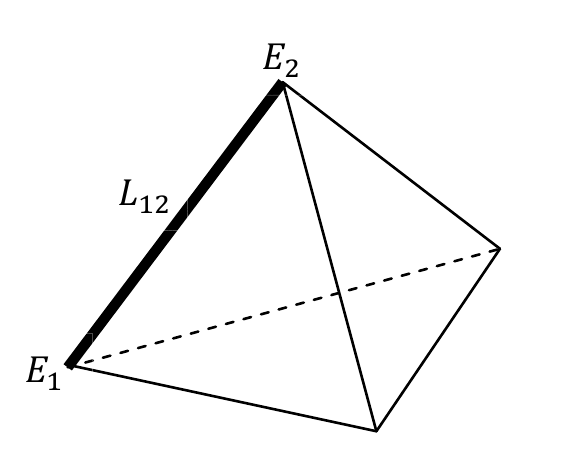}
\end{center}
$4$ elliptic fibrations coming from pencils of planes in ${\mathbb P}^3$ containing the line $R_{j}$, $j=1,2,3,4$.\\ 

\begin{center}
  \includegraphics{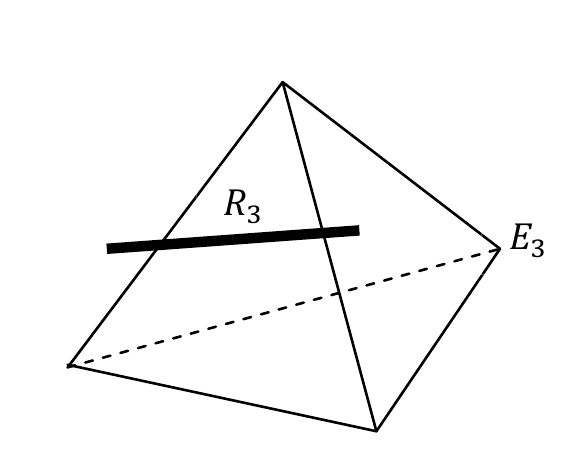}
\end{center}
One also has $4$ 'dual' elliptic fibrations coming from pencils of planes in ${\mathbb P}^3$ containing the line $E_{i}$ (on $X^{v}$), $i=1,2,3,4$. Another 4-tuple of pairs of 'dual' elliptic fibrations (corresponding to curves $E_i$ and $R_j$) can be constructed by using a different projective model of tetrahedral quartics - as quartics in ${\mathbb P}^3$ with $6$ coplanar nodes (see below) and considering pencils of planes in ${\mathbb P}^3$ containing a suitable line as above.\\ 

{\bf Lemma 0.8.} {\it Let $X$ be a general tetrahedral quartic.\\
Then each of the $6$ elliptic fibrations coming from the edges of the tetrahedron, has $2$ reducible fibers of type $I_4$ and $16$ fibers, which are nodal plane curves.\\
Each of the $4$ elliptic fibrations coming from residual lines (and their 'duals' coming from the nodes) has $1$ reducible fiber of type $I_6$ and $18$ fibers, which are nodal plane curves.\\
Each of the $4$ elliptic fibrations coming from residual lines and the $6$-nodal projective model of $X$ (and their 'duals' coming form the nodes) has $1$ reducible fiber of type $I_6$ and all the other singular fibers irreducible.}\\

{\it Proof:} Let us consider the elliptic fibration on $\tilde{X}$ coming from the line $L_{12}$. It is given by the complete linear system $\mid H-L_{12}-E_1-E_2 \mid$. It has two apparent reducible fibers (corresponding to the two faces of the tetrahedron adjacent at $L_{12}$): $L_{13}+L_{23}+E_3+R_4$ and $L_{14}+L_{24}+E_4+R_3$, each of which has type $I_4$ (i.e. is a cycle of $4$ $(-2)$-curves).\\

Similarly, the elliptic fibration on $\tilde{X}$ coming from the line $R_3$ is given by the complete linear system $\mid H - R_3 \mid$. It has one apparent reducible fiber (corresponding to the face of the tetrahedron containing $R_3$) $L_{12}+L_{14}+L_{24}+E_1+E_2+E_4$ of type $I_6$ (i.e. is a cycle of $6$ $(-2)$-curves).\\
 
For a general $X$ it follows from the explicit equation of the tetrahedral quartic and from the fact that a general pencil of plane cubics contains only smooth and nodal curves that all the irreducible singular fibers of these two types of elliptic fibrations on $X$ are nodal plane cubics.\\ 

It now follows from the comparison of Euler characteristics that there are $24-4-4=16$ of them in the first case and $24-6=18$ of them in the second case.\\

As for the elliptic fibrations coming from the $6$-nodal projective model of $X$, then we will see later that such $6$-nodal quartics contain exactly $4$ lines (if $X$ is general). So the elliptic fibrations coming from them (from the pencils of planes in ${\mathbb P}^3$ containing one of these lines) have exactly $1$ reducible fiber (corresponding to the plane in ${\mathbb P}^3$ containing all the $6$ nodes of the quartic) of type $I_6$. {\it QED}\\

{\bf Corollary 0.9.} {\it A general tetrahedral quartic $X$ contains exactly $10$ lines - the edges of the tetrahedron and the residual lines.}\\

{\it Proof:} Any line on $X$ would intersect a face of the tetrahedron at a point on one of the edges or at a point on a residual line. Hence it will be a component of a singular fiber of one of the elliptic fibrations considered above. But the components of the reducible fibers of those elliptic fibrations are nodes, edges of the tetrahedron and residual lines. {\it QED}\\

\subsection{A birational involution on the Hilbert square} 

Let $\pi\colon \tilde{X}\rightarrow X$ be the minimal resolution of singularities of a general tetrahedral quartic $X\subset {\mathbb P}^3$. Then any pair $P, P'$ of general points on $X$ determines a twisted cubic in ${\mathbb P}^3$ (which passes through $P, P'$ and the $4$ vertices of the tetrahedron), which in turn determines $2$ general points $Q, Q'$ on $X$ (such that this twisted cubic intersects $X$ at the $4$ vertices of the tetrahedron and the $4$ points $P, P', Q, Q'$). In other words, there exists a natural birational involution on the Hilbert square ${\tilde{X}}^{[2]}$ of the minimal resolution of singularities of any general tetrahedral quartic $X$.\\

\section{Neron-Severi lattice}

Let $U_1\subset {\mathbb C}^{12}$ be the (Zariski open) locus of points in the parameter space of tetrahedral quartics corresponding to tetrahedral quartics with exactly $4$ nodes (which in particular means that the residual lines do not intersect each other and nodes). We will call such tetrahedral quartics $4$-nodal.\\

Let $\pi\colon \tilde{X}\rightarrow X$ be the minimal resolution of singularities of a $4$-nodal tetrahedral quartic $X$.\\

Let us consider a lattice $M$ of rank $11$, which is defined by the incidence matrix of the curves $L_{12}, L_{13}, L_{14}, L_{23}, L_{24}, L_{34}, E_1, E_2, E_3, E_4, R_1$ on $\tilde{X}$:
$$
M=\left[
   \begin{array}{ccccccccccc}
-2 & 0 & 0 & 0 & 0 & 0 & 1 & 1 & 0 & 0 & 0\\
0 & -2 & 0 & 0 & 0 & 0 & 1 & 0 & 1 & 0 & 0\\
0 & 0 & -2 & 0 & 0 & 0 & 1 & 0 & 0 & 1 & 0\\
0 & 0 & 0 & -2 & 0 & 0 & 0 & 1 & 1 & 0 & 1\\
0 & 0 & 0 & 0 & -2 & 0 & 0 & 1 & 0 & 1 & 1\\
0 & 0 & 0 & 0 & 0 & -2 & 0 & 0 & 1 & 1 & 1\\
1 & 1 & 1 & 0 & 0 & 0 & -2 & 0 & 0 & 0 & 0\\
1 & 0 & 0 & 1 & 1 & 0 & 0 & -2 & 0 & 0 & 0\\
0 & 1 & 0 & 1 & 0 & 1 & 0 & 0 & -2 & 0 & 0\\
0 & 0 & 1 & 0 & 1 & 1 & 0 & 0 & 0 & -2 & 0\\
0 & 0 & 0 & 1 & 1 & 1 & 0 & 0 & 0 & 0 & -2
   \end{array}
 \right]
$$

One checks by a direct computation that $M$ is an even rank $11$ lattice with signature $(1,10)$ and discriminant $2^7$.\\

{\bf Lemma 1.1.} {\it If $X$ is a $4$-nodal tetrahedral quartic in ${\mathbb P}^3$ and $\pi\colon \tilde{X}\rightarrow X$ is its minimal resolution of singularities, then there is a primitive embedding of lattices $M\hookrightarrow NS(\tilde{X})$.}\\

{\it Proof:} If $l_{12}, l_{13}, l_{14}, l_{23}, l_{24}, l_{34}, e_1, e_2, e_3, e_4, r_1$ is the basis for $M$ as a ${\mathbb Z}$-module, in which its intersection product is given by the matrix above, then one identifies $l_{12}$ with $L_{12}$, $l_{13}$ with $L_{13}$, $l_{14}$ with $L_{14}$, $l_{23}$ with $L_{23}$, $l_{24}$ with $L_{24}$, $l_{34}$ with $L_{34}$, $e_1$ with $E_1$, $e_2$ with $E_2$, $e_3$ with $E_3$, $e_4$ with $E_4$ and $r_1$ with $R_1$. This gives an embedding of lattices $M\hookrightarrow NS(\tilde{X})$. One has to check that it is primitive.\\

Let $\hat{M}$ be the saturation of $M$ in $NS(\tilde{X})$. Then $[\hat{M}\colon M]^2=\frac{disc(M)}{disc(\hat{M})} \mid 2^7$. Hence $[\hat{M}\colon M]=2^{\alpha}$ and one needs to check that $\alpha =0$. For this we have to show that given a divisor class $D\in NS(\tilde{X})$ such that $2D\in M$, one has $D\in M$.\\

Let $2D=\sum_{1\leq i<j\leq 4}{a_{ij}\cdot L_{ij}}+\sum_{i=1}^{4}{b_i\cdot E_i}+c\cdot R_1$, where $a_{ij},b_i,c\in\{ 0,1 \}$. We will be working in the group $NS(\tilde{X})/ 2\cdot NS(\tilde{X})$.\\

If all the $a_{ij}$ are zero, then all $b_i$ are zero and $c=0$ as well by a theorem of Nikulin saying that $k$ smooth irreducible disjoint rational curves on a $K3$ surface can represent an even divisor class in the Neron-Severi group only if $8\mid k$ \cite{NikulinKummer}.\\

So, we may assume that $a_{14}=1$. Since the intersection of $2D$ with $E_1$ is even, we must have that $a_{13}+a_{12}=1$. We may assume that $a_{13}=1$ and $a_{12}=0$. Since the intersection of $2D$ with $R_2$ is even, we must have $a_{34}=0$. Similarly, we consider intersections of $2D$ with $E_3$ and $E_4$ and conclude that $a_{23}=a_{24}=1$.\\

So, we get: $2D=L_{13}+L_{14}+L_{23}+L_{24}+\sum_{i=1}^{4}{b_i\cdot E_i}+c\cdot R_1$.\\

Intersecting $2D$ with $L_{12}$ gives: $b_1+b_2=0\; mod\; 2$, i.e. $b_1=b_2$. Similarly, intersecting $2D$ with $L_{13}$ and $L_{14}$ we get: $b_1=b_2=b_3=b_4=b$. Since $0=2D\cdot L_{34}=b_3+b_4+c=c \; mod \; 2$, we get $c=0$.\\

If $b=0$, then we arrive at contradiction by Nikulin's theorem, quoted above. Hence $b=1$.\\

So, we get: $2D=L_{13}+L_{14}+L_{23}+L_{24}+E_1+E_2+E_3+E_4$ in $NS(\tilde{X})$. Let us check that this is impossible.\\ 

We notice that $D^2=0$. Hence $D$ is effective (if $-D$ were effective, then $-2D$ would be also effective, which is not the case).\\

If $E_i$ were in the fixed locus of $\mid D \mid$, then for any $j\neq i$ we would have $1\leq (2D)\cdot L_{ij}=-2+2=0$ (we can arrange that $(ij)\neq (12)$ and $(ij)\neq (34)$), which is impossible. Hence if $\Gamma$ is a nonsingular rational curve in the base locus of $\mid D \mid$, then its image in ${\mathbb P}^3$ ${\pi}(\Gamma)$ would be a curve. Since $D\cdot H=\frac{1}{2}\cdot (2D\cdot H)=\frac{4}{2}=2$, the base locus of $\mid D \mid$ consists of at most two $(-2)$-curves.\\

If $D={\Gamma}_1+{\Gamma}_2$, where ${\Gamma}_1$, ${\Gamma}_2$ are $(-2)$-curves and ${\pi}({\Gamma}_1)$ and ${\pi}({\Gamma}_2)$ are lines in ${\mathbb P}^3$, then ${\Gamma}_1 \cdot {\Gamma}_2=2$, which is impossible for lines in ${\mathbb P}^3$. So, either $\mid D\mid$ is base point free, or the base locus of $\mid D\mid$ consists of one $(-2)$-curve. In the latter case, $D=E+\Gamma$, where $dim \mid E\mid\geq 1$ and $\mid E\mid$ is base point free. Since ${\pi}(E)\cdot H={\pi}(D)\cdot H-{\pi}(\Gamma)\cdot H=2-1=1$, ${\pi}(E)$ should be a line on $X$, i.e. $dim \mid E\mid=0$. This is a contradiction.\\

Hence $\mid D\mid$ is base point free. So, by \cite{SD} $D=k\cdot E$, where $k\geq 1$ and $E$ is an elliptic pencil on $\tilde{X}$. Since $E$ is nef and $0=(2D)\cdot E=L_{13}\cdot E+L_{14}\cdot E+L_{23}\cdot E+L_{24}\cdot E+E_1\cdot E+E_2\cdot E+E_3\cdot E+E_4\cdot E$ we get: $0=L_{13}\cdot E=L_{14}\cdot E=L_{23}\cdot E=L_{24}\cdot E=E_1\cdot E=E_2\cdot E=E_3\cdot E=E_4\cdot E$, i.e. the curves $L_{13}, L_{14}, L_{23}, L_{24}, E_1, E_2, E_3, E_4$ lie in the fibers of the elliptic fibration ${\phi}_{\mid E\mid}\colon \tilde{X}\rightarrow {\mathbb P}^1$ defined by $\mid E\mid$. Since the union of supports of these curves is connected, they should lie in the same fiber of ${\phi}_{\mid E\mid}$, i.e. $E=L_{13}+L_{14}+L_{23}+L_{24}+E_1+E_2+E_3+E_4+D'$ for some effective divisor $D'$. Hence $E\geq 2D=2k\cdot E$, which is impossible.\\

So, we conclude that $D\in M$, and hence $M\subset NS(\tilde{X})$ is a primitive sublattice. {\it QED}\\

In order to compute the Neron-Severi lattice of a general tetrahedral quartic, we will use the theory of lattice polarized $K3$ surfaces and their moduli from \cite{Dolgachev}. Note, that a general tetrahedral quartic $X$ (or its minimal resolution of singularities $\tilde{X}$, to be more precise) is $M$-polarized. Let us also note that if $\Lambda=E_8(-1)^{{\oplus} 2}{\oplus} H^{{\oplus} 3}$ denotes the cohomology lattice $H^2(Y,\mathbb Z)$ of a $K3$ surface $Y$, then there exists a unique primitive embedding of lattices $M\hookrightarrow \Lambda$. This is a consequence of a criterion due to Nikulin \cite{Nikulin}. Indeed, one computes directly that the discriminant-group of $M$ is $A_M=M^{*}/M\cong ({\mathbb Z}/8{\mathbb Z})\oplus ({\mathbb Z}/4{\mathbb Z})^{{\oplus}2}$ and it is generated by $l(A_M)=l((A_M)_2)=3$ elements, while $rk(M)=11$.\\ 

{\bf Lemma 1.2.} {\it If $X$ is a $4$-nodal tetrahedral quartic, then divisors $A=L_{12}+L_{13}+L_{14}+L_{23}+L_{24}+L_{34}+E_1+E_2+E_3+E_4+R_1+R_2+R_3+R_4$ and $A_0=3H-E_1-E_2-E_3-E_4$ in $M\subset NS(\tilde{X})$ are very ample.}\\

{\it Proof:} We apply Nakai-Moishezon criterion. Let us consider the case of $A$. $A^2=20>0$, $A\cdot E_i=A\cdot R_j=1>0$, $A\cdot L_{ij}=2>0$. If $C\subset \tilde{X}$ is an irreducible curve different from lines and nodes, then $H\cdot C=deg {\pi}(C)>0$. Since $A=H+D$, where $D\geq 0$ and $Supp(D)\subset {\cup}_{i,j}L_{ij}\cup {\cup}_{i}E_{i}\cup {\cup}_{j}R_{j}$, we have $D\cdot C\geq 0$, and so $A\cdot C=H\cdot C+D\cdot C>0$. Hence $A$ is ample. Similarly, one checks that $A_0$ is ample.\\ 

Very ampleness of $A$ and $A_0$ follows from the criteria of Saint-Donat (see \cite{SD}, \cite{Mori}). Indeed, for any ample divisor $A'$ such that $\mid A'-H \mid \neq \emptyset$, $(A')^2\geq 4$, $(A')^2\neq 8$ (both $A$ and $A_0$ satisfy these conditions) one sees that the conditions of Theorem 5 in \cite{Mori} are satisfied. So, $A'$ is very ample. {\it QED}\\

{\bf Proposition 1.3.} {\it Let $X$ be a very general tetrahedral quartic in ${\mathbb P}^3$ (i.e. $X$ corresponds to a point in ${\mathbb C}^{12}$ outside of a countable union of algebraic hypersurfaces, and in particular we may assume that $X$ has exactly $4$ nodes and exactly $10$ lines) and $\pi\colon \tilde{X}\rightarrow X$ be its minimal resolution of singularities. Then the primitive embedding $M\hookrightarrow NS(\tilde{X})$ is an isomorphism.}\\

{\it Proof:} Let $D_M/{\Gamma}_M$ be the coarse moduli space of $M$-polarized $K3$ surfaces \cite{Dolgachev} and $\tau\colon U_1\rightarrow D_M/{\Gamma}_M$ be the morphism corresponding to the natural family of tetrahedral quartics $\mathcal X\rightarrow U_1$ (or their simultaneous minimal resolution of singularities $\tilde{\mathcal X}\rightarrow U_1$, where $\tilde{\mathcal X}$ is obtained form $\mathcal X$ by blowing-up the loci of nodes).\\

$M$-polarized $K3$ surfaces $Y$ with $NS(Y)=M$ include those, which correspond to the analytically open subset $W\subset D_M/{\Gamma}_M$, which is the complement of a countable union of algebraic hypersurfaces in $D_M/{\Gamma}_M$. So, it is sufficient to check that $W$ lies in the image of $\tau$.\\

Let $Y$ be a $K3$ surface with $NS(Y)\cong M$ and $X$ be a tetrahedral quartic with $4$ nodes and $10$ lines corresponding to a point in $U_1\subset {\mathbb C}^{12}$, $\pi\colon \tilde{X}\rightarrow X$ be its minimal resolution of singularities. We denote by $L_{ij}, E_i, R_j, H, A, A_0$ the elements of $M$ and the corresponding divisor classes on $Y$ under an isomorphism $NS(Y)\cong M$, and by $\tilde{L_{ij}}, \tilde{E_i}, \tilde{R_j}, \tilde{H}, \tilde{A}, \tilde{A_0}$ the corresponding divisor classes on $\tilde{X}$ under the primitive embedding $M\hookrightarrow NS(\tilde{X})$.\\

Since $A^2>0$ and the nef cone of a $K3$ surface is the fundamental domain for the action of the group generated by reflections in $(-2)$-classes on the positive cone of the $K3$ surface, we may assume (at the cost of changing the isomorphism $NS(Y)\cong M$) that $A$ is a nef divisor on $Y$. In fact, it is an ample divisor on $Y$, since if there existed a $(-2)$-class $\delta \in M\cong NS(Y)$ such that $A\cdot \delta =0$, then there would be also a $(-2)$-class $\tilde{\delta}\in NS(\tilde{X})\hookleftarrow M$ such that $\tilde{A}\cdot \tilde{\delta} =0$, which is impossible, since $\tilde{A}$ is ample. So, $A$ is an ample divisor on $Y$.\\

Since $E_i\cdot A=R_j\cdot A=1$, each divisor class $E_i$ and $R_j$ is represented by a smooth irreducible rational curve on $Y$. Since $L_{ij}\cdot A=2$, each divisor class $L_{ij}$ is represented either by a smooth irreducible rational curve on $Y$, or is a sum of two irreducible divisors: $L_{ij}=D_i+D_j$. The latter is impossible, since $NS(Y)\cong M\hookrightarrow NS(\tilde{X})$. Indeed, if $\tilde{L_{ij}}=\tilde{D_i}+\tilde{D_j}$, then $\tilde{D_i}$ and $\tilde{D_j}$ are represented by $(-2)$-curves on $\tilde{X}$ and $\tilde{A}\cdot \tilde{D_i}=\tilde{A}\cdot \tilde{D_j}=1$. Hence ${\pi}(\tilde{D_i})$ and ${\pi}(\tilde{D_j})$ are either nodes or the lines on $X$, i.e. one of the divisors $\tilde{L_{12}}, \tilde{L_{13}}, \tilde{L_{14}}, \tilde{L_{23}}, \tilde{L_{24}}, \tilde{L_{34}}, \tilde{E_1}, \tilde{E_2}, \tilde{E_3}, \tilde{E_4}, \tilde{R_1}, \tilde{R_2}, \tilde{R_3}, \tilde{R_4}$. This is impossible, as one sees immediately from considering intersection products. So, all the divisor classes $L_{ij}, E_{i}, R_{j}$ are represented by smooth irreducible rational curves on $Y$.\\

Let us consider the nef divisor $H=L_{12}+L_{13}+L_{23}+E_1+E_2+E_3+R_4$. Since $L_{ij}, E_{i}, R_{j}$ are $14$ distinct irreducible curves on $Y$, and $H=L_{13}+L_{14}+L_{34}+E_1+E_3+E_4+R_2=$ $L_{12}+L_{14}+L_{24}+E_1+E_2+E_4+R_3=$ $L_{23}+L_{24}+L_{34}+E_2+E_3+E_4+R_1$ we see that $\mid H\mid$ has no fixed components. Hence it is base point free \cite{SD}. In particular, since $H^2=4$, we have $h^1(Y,H)=0$, $dim \mid H \mid=3$ and $\mid H \mid$ defines a morphism ${\phi}_{\mid H \mid}\colon Y\rightarrow {\mathbb P}^3$, which is a birational morphism of $Y$ onto a quartic surface in ${\mathbb P}^3$ and is an isomorphism except that it contracts $E_1, E_2, E_3, E_4$ to the $4$ nodes ${\phi}_{\mid H \mid}(E_1)$, ${\phi}_{\mid H \mid}(E_2)$, ${\phi}_{\mid H \mid}(E_3)$, ${\phi}_{\mid H \mid}(E_4)$ on ${\phi}_{\mid H \mid}(Y)$. These nodes are connected by the $6$ lines ${\phi}_{\mid H \mid}(L_{ij})$, $1\leq i<j\leq 4$, which form a tetrahedron in ${\mathbb P}^3$. Hence ${\phi}_{\mid H \mid}(Y)$ is a tetrahedral quartic in ${\mathbb P}^3$ (corresponding to a point in $U_1\subset {\mathbb C}^{12}$) and ${\phi}_{\mid H \mid}\colon Y\rightarrow {\phi}_{\mid H \mid}(Y)$ is its minimal resolution of singularities.\\

This means that any point of $W\subset D_M/{\Gamma}_M$ lies in the image of $\tau\colon U_1\rightarrow D_M/{\Gamma}_M$, as we wanted. {\it QED}\\

{\bf Remark.} One easily sees that ${\phi}_{\mid H \mid}(Y)$ contains exactly $10$ lines. Indeed, if $C\subset Y$ is a strict transform of a line on ${\phi}_{\mid H \mid}(Y)$, then $H\cdot C=1$, and so the corresponding divisor class $\tilde{C}\in NS(\tilde{X})\hookleftarrow M$ is represented by a line on $X$. Hence $\tilde{C}$ is one of the lines $\tilde{L_{ij}}, \tilde{R_{k}}$. So, $C$ is one of the lines $L_{ij}, R_k$.\\

Let $W_1\subset D_M/{\Gamma}_M$ be the subset corresponding to $M$-polarized $K3$ surfaces $Y$ such that $A\in M$ and $A_0\in M$ are ample divisors on $Y$. By Lemma 1.2, $W_1$ contains the image of the morphism $\tau \colon U_1\rightarrow D_M/{\Gamma}_M$ introduced above. We will see that in fact $Im(\tau)=W_1$ is a Zariski open dense subset of $D_M/{\Gamma}_M$ - the coarse moduli space of $4$-nodal tetrahedral quartics.\\

{\bf Lemma 1.4.} {\it Let $Y$ be an $M$-polarized $K3$ surface such that $A\in M\hookrightarrow NS(Y)$ is an ample divisor. Then $A_0\in M\hookrightarrow NS(Y)$ is nef.}\\

{\it Proof:} Let us recall that $A=L_{12}+L_{13}+L_{14}+L_{23}+L_{24}+L_{34}+E_1+E_2+E_3+E_4+R_1+R_2+R_3+R_4$ and $A_0=3H-E_1-E_2-E_3-E_4=$ $2E_4+E_1+E_2+E_3+R_1+R_2+R_3+2(L_{14}+L_{24}+L_{34})+L_{12}+L_{23}+L_{13}$\\  

Since $A\cdot E_i=A\cdot R_j=1$ for any $i$ and $j$ and $A$ is ample by assumption, $E_1, E_2, E_3, E_4$, $R_1, R_2, R_3, R_4$ are represented by $8$ disjoint smooth irreducible rational curves on $Y$. Since $A\cdot L_{ij}=2$, $L_{ij}$ is either irreducible (and in this case is represented by a smooth irreducible rational curve different from $E_1, E_2, E_3, E_4$, $R_1, R_2, R_3, R_4$), or $L_{ij}=D_i+D_j$, where $D_i$ and $D_j$ are irreducible and $A\cdot D_i=A\cdot D_j=1$. Let us check that $D_i$ and $D_j$ are represented by smooth irreducible rational curves and $Supp(L_{ij})\cap Supp(L_{i'j'})=\emptyset$, if $(ij)\neq(i'j')$.\\

If $(D_i)^2\geq 0$, then $dim \mid D_i \mid \geq 1$. Hence since $D_i$ is irreducible, $\mid D_i \mid$ has no fixed components and so is base point free. Let us consider the corresponding morphism ${\phi}_{\mid D_i \mid}\colon Y\rightarrow {\mathbb P}^N$, where $N=dim \mid D_i \mid \geq 1$. Note that $D_i$ is nef. Since $1=A\cdot D_i$, we have equalities $L_{ab}\cdot D_i=0$, $E_c\cdot D_i=0$, $R_d\cdot D_i=0$ for all the divisors $L_{ab}$, $E_c$, $R_d$, except for exactly one. Since the union of supports of remaining $13$ divisors is connected, they all should lie in the same fiber of ${\phi}_{\mid D_i \mid}$. In particular, say, $E_1, E_2, E_3, R_1, R_2, R_3$ lie in one fiber of ${\phi}_{\mid D_i \mid}$, i.e. $D_i \geq E_1+E_2+E_3+R_1+R_2+R_3$. Since $A$ is ample and $A\cdot E_i=A\cdot R_j=1$, this implies that $1=A\cdot D_i\geq 6$. Contradiction.\\

Hence $(D_i)^2=(D_j)^2=-2$, i.e. $D_i$ and $D_j$ are represented by smooth irreducible rational curves on $Y$. Since $(L_{ij})^2=-2$, we get: $D_i\cdot D_j=1$.\\

\begin{center}
  \includegraphics{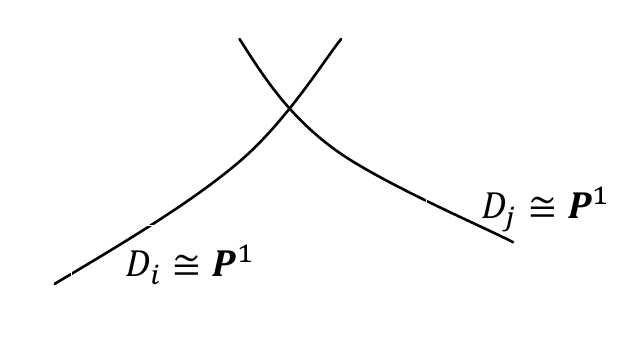}
\end{center}
Note that irreducible components of $L_{ij}$ are not among the curves $E_1, E_2, E_3, E_4$, $R_1, R_2, R_3, R_4$, because $D_i\cdot L_{ij}=-1$.\\

Suppose $D=D_{i_1}=D_{i_2}$ is a common irreducible component of $L_{i_1j_1}$ and $L_{i_2j_2}$. Let $L_{i_3j_3}$ be the third divisor, which shares an irreducible component with $L_{i_1j_1}+L_{i_2j_2}$.\\

\begin{center}
  \includegraphics{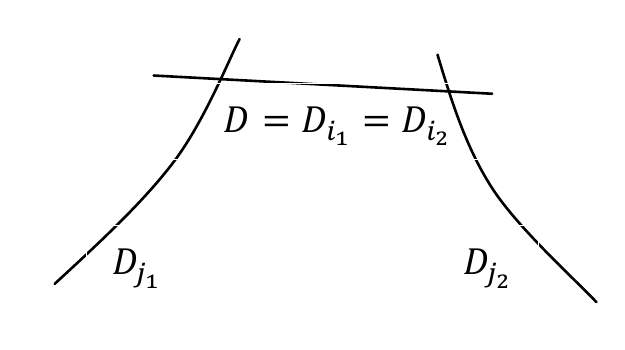}
\end{center}
Then this irreducible component must be $D=D_{i_1}=D_{i_2}=D_{i_3}$, because otherwise (if, say, $D_{j_2}=D_{j_3}$), $0=L_{i_3j_3}\cdot L_{i_1j_1}\geq 1$, which is impossible.\\

\begin{center}
  \includegraphics{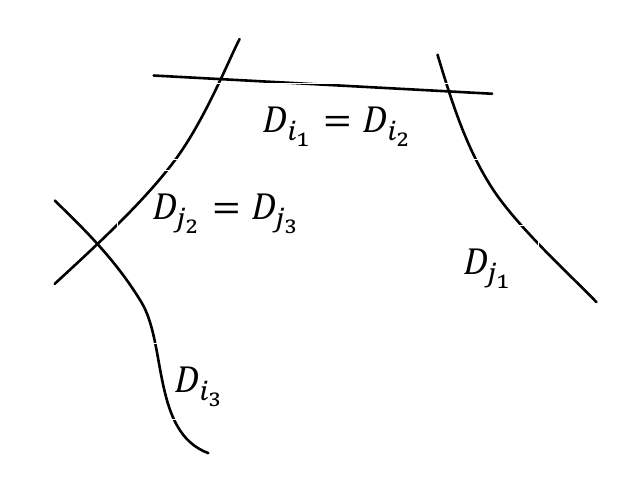}
\end{center} 
Note that since $(L_{ij}\cdot L_{i'j'})=0$, either $L_{ij}$ and $L_{i'j'}$ share an irreducible component, or their supports are disjoint. Hence a connected component of the support of the divisor $\sum_{1\leq i<j\leq 4}{L_{ij}}$ should have the form $L_{i_1j_1}+\dots + L_{i_kj_k}$ and should look as on the picture:\\

\begin{center}
  \includegraphics{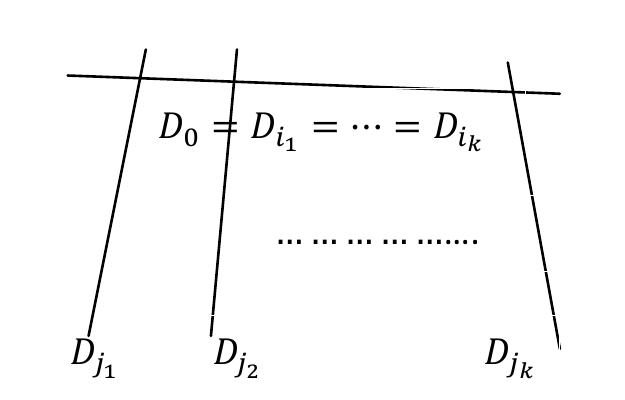}
\end{center} 
We have: $1=A\cdot D_0=-k+\sum_{i}(E_i\cdot D_0)+\sum_{j}(R_j\cdot D_0)$. Hence the set $\{ i\; \mid \; E_i\cdot D_{i_{\alpha}}=1 \} \cup \{ 4+j\; \mid \; R_j\cdot D_{i_{\alpha}}=1 \}\subset \mathbb N$ does not depend on $\alpha =1,...,k$. Hence it contains two elements (corresponding to $E_i$ and $R_j$ - look at $L_{12}$ and $L_{13}$, for example). In particular, $k=1$. Hence if $(ij)\neq(i'j')$, then $Supp(L_{ij})\cap Supp(L_{i'j'})=\emptyset$.\\

Since $A_0$ is effective, in order to check that $A_0$ is nef it is sufficient to consider irreducible curves $C\subset Supp(A_0)$. Since $A_0\cdot E_i=2\geq 0$, $A_0\cdot R_j=3\geq 0$, $A_0\cdot L_{ij}=1\geq 0$, it is enough to check that $A_0.\cdot D_i\geq 0$, when $L_{ij}=D_i+D_j$. Since $D_i\cdot L_{ij}=-1$ and $D_i\cdot L_{i'j'}=0$, if $(i'j')\neq (ij)$, we get (we may assume that $L_{ij}=L_{12}$ without loss of generality): $A_0\cdot D_i\geq 1-1=0$.\\

Note that $D_i\cdot (E_1+E_2+E_3+R_1+R_2+R_3)\geq 1$, since otherwise $D_i\cdot E_1=D_i\cdot E_2=D_i\cdot E_3=D_i\cdot R_1=D_i\cdot R_2=D_i\cdot R_3=D_i\cdot L_{ab}=0$ for $(ab)\neq(1,2)$, and so $1=D_i\cdot A=L_{12}\cdot D_i+R_4\cdot D_i=-1+1=0$, which is impossible.\\

So, $A_0$ is nef. {\it QED}\\

{\bf Corollary 1.5.} {\it The subset $W_1\subset D_M/{\Gamma}_M$ is Zariski open.}\\

{\it Proof:} Since $A^2=20>0$ and the primitive embedding $M\hookrightarrow \Lambda$ is unique, we may assume that $A$ is nef. Then the condition that $A$ is ample on an $M$-polarized $K3$ surface $Y$ may be rephrased as follows: there is no $\delta \in NS(Y)$ such that $(\delta)^2=-2$ and $A\cdot\delta=0$. If we impose an analogous condition for $A_0$, this will ensure (thanks to the Lemma 1.4) that $A_0$ is ample as well.\\

Let us consider lattice $N=M^{\bot}_{\Lambda}$ of signature $(2,9)$. Then $D_M=\{ \omega \in {\mathbb P}(N_{\mathbb C})\; \mid \; {\omega}\cdot{\omega}=0, \; {\omega}\cdot \bar{\omega}>0 \}$. Let $N=[\Lambda \colon N\oplus M]$. Then any $\delta \in \Lambda$ such that $(\delta)^2=-2$ and ${\delta}\cdot A=0$ or ${\delta}\cdot A_0=0$ can be written as $\delta = \frac{1}{a}\cdot (m+n)$, where $1\leq a\leq N$, $n\in N$, $n\neq 0$, $m\in M$, $A\cdot m=0$ or $A_0\cdot m=0$.\\

Hence $-2a^2=m^2+n^2$. Since $sign(M)=(1,10)$ and $A^2>0$ (respectively $(A_0)^2>0$) we have $m^2\leq 0$ and $n^2=-2a^2-m^2\geq -2a^2\geq -2N^2$. If $n^2\geq 0$, then $\pm n$ is an effective divisor, and so it can not lie on an ample $M$-polarized $K3$ surface $Y$.\\

Let us consider elements $n\in N$ and $m\in M$ such that $-1\geq n^2\geq -2N^2$, $m\cdot A=0$ or $m\cdot A_0=0$, $m^2\geq 1-2N^2$ and there exists $\delta \in \Lambda$ and $a\geq 1$, $a\leq N$ such that $(\delta)^2=-2$ and $a\cdot \delta=m+n$. There are finitely many such elements $m\in M$ and $n\in N$ upto the action of $O(N)$ by the Remark on page 2607 in \cite{Dolgachev}. For any $n\in N$ let us denote by $H_n$ the hyperplane $H_n=\{ \omega \in D_M  \;  \mid \;  {\omega}\cdot n=0 \}$ and the corresponding hypersurface in $D_M/{\Gamma}_M$. We conclude (as in Remark in \cite{Dolgachev} quoted above) that there are finitely many algebraic hypersurfaces $H_{n_1},\dots , H_{n_k}\subset D_M/{\Gamma}_M$ such that their complement is exactly the subset $W_1\subset D_M/{\Gamma}_M$. {\it QED}\\

{\bf Theorem 1.6.} {\it $Im(\tau)=W_1$, i.e. if $Y$ is an $M$-polarized $K3$ surface such that $A\in M$ and $A_0\in M$ correspond to ample divisor classes, then $Y$ is (the minimal resolution of singularities of ) a $4$-nodal tetrahedral quartic in ${\mathbb P}^3$.}\\

{\it Proof:} Since $A$ and $A_0$ are ample and $A\cdot E_i=A\cdot R_j=A_0\cdot L_{ij}=1$, $E_1, E_2, E_3, E_4$, $R_1, R_2, R_3, R_4$, $L_{12}, L_{13}, L_{14}, L_{23}, L_{24}, L_{34}$ are $14$ distinct smooth irreducible rational curves on $Y$.\\

Hence $H=L_{12}+L_{13}+L_{23}+E_1+E_2+E_3+R_4=$ $L_{13}+L_{14}+L_{34}+E_1+E_3+E_4+R_2=$ $L_{12}+L_{24}+L_{14}+E_1+E_2+E_4+R_3=$ $L_{24}+L_{23}+L_{34}+E_2+E_3+E_4+R_1$ is a nef divisor with base point free complete linear system and such that $H^2=4$. So, its complete linear system determines a morphism ${\phi}_{\mid H \mid}\colon Y\rightarrow {\mathbb P}^3$, which represents $Y$ as the minimal resolution of singularities of a tetrahedral quartic in ${\mathbb P}^3$ ($E_1, E_2, E_3, E_4$ get contracted into the $4$ nodes and $L_{ij}$ map to the lines connecting them, i.e. to the edges of a tetrahedron).\\

Note that $\mid H \mid$ is not hyperelliptic (in the terminology of \cite{SD}) by Theorem 5.2 in \cite{SD}. Indeed, otherwise by that theorem there would exist a smooth irreducible curve $E$ of genus $1$ on $Y$ such that $E\cdot H=2$. This leads to an immediate contradiction, if one considers the elliptic fibration ${\phi}_{\mid E \mid}\colon Y\rightarrow {\mathbb P}^1$ associated with $E$ and the intersection products of $E$ with divisors $L_{ij}$, $E_i$, $R_j$, $H$. {\it QED}\\

{\bf Remark.} (1) We saw earlier that if $\pi\colon \tilde{X}\rightarrow X$ is the minimal resolution of singularities of a $4$-nodal tetrahedral quartic, then $A, A_0\in M\subset NS(\tilde{X})$ are very ample.\\

The closed immersion ${\phi}_{\mid A \mid}\colon \tilde{X}\rightarrow {\mathbb P}^{11}$ represents $\tilde{X}$ as a nondegenerate surface of degree $20$ in ${\mathbb P}^{11}$ containing $8$ disjoint lines $E_1, E_2, E_3, E_4$, $R_1, R_2, R_3, R_4$ and $6$ disjoint smooth conics $L_{12}, L_{13}, L_{14}, L_{23}, L_{24}, L_{34}$ such that each conic intersects $4$ lines and each line intersects $3$ conics. It 'squeezes' the tetrahedron in a sense that there is a hyperplane section of ${{\phi}_{\mid A \mid}}(\tilde{X})$, which coincides with $(Tetrahedron)\cap \tilde{X}$.\\

The closed immersion ${\phi}_{\mid A_0 \mid}\colon \tilde{X}\rightarrow {\mathbb P}^{15}$ represents $\tilde{X}$ as a nondegenerate surface of degree $28$ in ${\mathbb P}^{15}$ containing $6$ disjoint lines $L_{12}, L_{13}, L_{14}, L_{23}, L_{24}, L_{34}$, $4$ disjoint smooth conics $E_1, E_2, E_3, E_4$ and $4$ disjoint twisted cubics $R_1, R_2, R_3, R_4$ such that twisted cubics and conics do not intersect, each line intersects $2$ conics and $2$ twisted cubics, each conic intersects $3$ lines. Each twisted cubic intersects $3$ lines.\\

(2) Theorem 1.6. describes Neron-Severi lattices of (the minimal resolutions of singularities of) $4$-nodal (in particular, general) tetrahedral quartics. They are exactly those lattices $L$, which are Neron-Severi lattices of $K3$ surfaces and where lattice $M$ can be primitively embedded such that $A$ and $A_0$ are not orthogonal to elements of $L$ with square $-2$.\\

A $K3$ surface $Y$ is (the minimal resolution of singularities of ) a $4$-nodal tetrahedral quartic, if and only if it can be $M$-polarized (i.e. a primitive embedding $M\hookrightarrow NS(Y)$ can be chosen) in such a way that $A$ and $A_0$ become ample divisors. This characterization of $4$-nodal tetrahedral quartics is in the spirit of the characterization of Kummer surfaces via lattice polarizations (see \cite{NikulinKummer}).\\ 

\section{Mori-Mukai uniqueness} 

There are several ways to distinguish tetrahedral quartics (among each other and among other $K3$ surfaces) by looking at the curves naturally associated to them.\\

Let $\pi\colon \tilde{X}\rightarrow X$ be the minimal resolution of singularities of a $4$-nodal tetrahedral quartic $X\subset {\mathbb P}^3$. Let $C=A+L_{12}$. Then $C$ is nef, $C^2=22$ and the complete linear system $\mid C\mid$ is base point free, because otherwise $L_{12}$ would be a fixed component of $\mid C\mid$, which is impossible by criteria of Saint-Donat (see \cite{SD} or \cite{Mori}). Hence by another theorem of Saint-Donat \cite{SD} (and Bertini's theorem) general elements $Z_1\in \mid A\mid$ and $Z_2\in \mid C\mid$ are smooth irreducible curves on $\tilde{X}$ of genera $g(Z_1)=1+\frac{A^2}{2}=11$ and $g(Z_2)=1+\frac{C^2}{2}=12$, which intersect transversely in $A\cdot C=22$ points. In particular, the arithmetic genus of the nodal curve $Z_1\cup Z_2$ is $p_a(Z_1\cup Z_2)=44$.\\

{\bf Proposition 2.1.} {\it Let $\pi\colon \tilde{X}\rightarrow X$ be the minimal resolution of singularities of a $4$-nodal tetrahedral quartic $X\subset {\mathbb P}^3$. Let $C=A+L_{12}$ and $Z_1\in \mid A\mid$, $Z_2\in \mid C\mid$ be general elements. Then the nodal curve $Z_1\cup Z_2$ on $\tilde{X}$ distinguishes $\tilde{X}$ among all $K3$ surfaces in the following sense:}
\begin{itemize}
\item[(*)] {\it $\tilde{X}$ is the unique (upto isomorphism) $K3$ surface containing the curve $Z_1\cup Z_2$ such that the class of its divisor is ample and not divisible in $NS(\tilde{X})$}
\end{itemize}

{\bf Remark.} In particular, $\tilde{X}$ is the unique (upto isomorphism) $K3$ surface, which is the minimal resolution of singularities of a $4$-nodal tetrahedral quartic and contains the curve $Z_1\cup Z_2$ in the linear system of the divisor $2A+L_{12}$.\\ 

{\it Proof:} We will use Corollary 1.6 from \cite{MoriMukai}.\\

$A+C$ is not divisible in $NS(\tilde{X})$, because $(A+C)\cdot L_{12}=2$ and $(A+C)\cdot E_1=3$ are coprime.\\

Let us check that conditions (1.0) of \cite{MoriMukai} are satisfied. Recall, that $A$ is very ample and so its complete linear system defines a closed immersion $\tilde{X}\hookrightarrow {\mathbb P}^m$ with $m=11$ and hyperplane section class $A$. It is also immediate by Nakai-Moishezon criterion that $A+C$ is ample.\\

We need to check the following:
\begin{itemize}
\item[(i)] $\tilde{X}$ in ${\mathbb P}^m$ is cut out by quadrics
\item[(ii)] restriction homomorphism $H^0(\tilde{X}, A) \rightarrow H^0(C,A{\mid}_C)$ is an isomorphism
\item[(iii)] $deg_C(A{\mid}_C)=A\cdot C\geq m+1$
\end{itemize}

Condition $(iii)$ is immediate, because $m=1+\frac{A^2}{2}=11$ and $A\cdot C=22$.\\

Condition $(i)$ follows from Theorem 7.2 of \cite{SD}. Indeed, if there were a smooth irreducible curve $E$ on $\tilde{X}$ such that $E^2=0$ and $E\cdot A=3$, then we would have: $E\cdot H=3$ (since otherwise $E$ would have to be rational, which is impossible) and $0=E\cdot(A-H)=E\cdot E_3+E\cdot R_1+E\cdot R_2+E\cdot R_4+E\cdot L_{23}+E\cdot L_{34}+E\cdot L_{13}$. Since $E$ is nef, this implies that $E\cdot E_3=E\cdot R_1=E\cdot R_2=E\cdot R_4=E\cdot L_{23}=E\cdot L_{34}=E\cdot L_{13}=0$. Since the complete linear system $\mid E\mid$ is base point free and $Supp(E_3+R_1+R_2+R_4+L_{23}+L_{34}+L_{13})$ is connected, this implies that it lies in a fiber of ${\phi}_{\mid E \mid}\colon \tilde{X}\rightarrow {\mathbb P}^N$, where $N=dim \mid E \mid \geq 1$. Hence $E\geq E_3+R_1+R_2+R_4+L_{23}+L_{34}+L_{13}$, and so $3=E\cdot H\geq 6$. Contradiction.\\  

Let us check condition $(ii)$. Since $\mid A-C\mid = \mid -L_{12} \mid = \emptyset$, the image of $Z_2$ under the embedding $\tilde{X}\hookrightarrow {\mathbb P}^m$ is a nondegenerate curve. So, it is sufficient to check that $h^0(\tilde{X},A)=h^0(Z_2,A{\mid}_{Z_2})$. By \cite{SD}, Proposition 2.6 $h^1(\tilde{X}, A)=0$. Hence by Riemann-Roch theorem $h^0(\tilde{X},A)=2+\frac{A^2}{2}=12$.\\

By adjunction formula $K_{Z_2} = C{\mid}_{Z_2}$. Hence by Serre duality $h^1(Z_2,A{\mid}_{Z_2})=$ $h^0(Z_2,(C-A){\mid}_{Z_2})=$ $h^0(Z_2,{\mathcal O}_{\tilde{X}}(L_{12}){\mid}_{Z_2})$. Since $L_{12}\cdot C=0$ and $\mid C \mid$ is base point free, this implies that $Z_2$ does not intersect $L_{12}$,i.e. ${\mathcal O}_{\tilde{X}}(L_{12}){\mid}_{Z_2}\cong {\mathcal O}_{Z_2}$ and so $h^1(Z_2, A{\mid}_{Z_2})=1$. Hence by Riemann-Roch theorem $h^0(Z_2, A{\mid}_{Z_2})=A\cdot C-g(Z_2)+1+h^1(Z_2, A{\mid}_{Z_2})=22-12+1+1=12$, i.e. $h^0(Z_2,A{\mid}_{Z_2})=h^0(\tilde{X},A)$. {\it QED}\\

In fact, Corollary 1.6 of \cite{MoriMukai} says more:\\ 

{\bf Corollary 2.2.} {\it Let $Z\subset \tilde{X}$ be a Deligne-Mumford stable curve such that $p_a(Z)=44$ and the corresponding divisor class is ample and not divisible in $NS(\tilde{X})$.\\

Then any isomorphism $\gamma \colon Z \rightarrow Z_1\cup Z_2$ (of abstract nodal curves) can be extended to an automorphism $\Gamma \colon \tilde{X} \rightarrow \tilde{X}$ such that ${\Gamma}{\mid}_{Z}=\gamma$ and ${\Gamma}(Z)=Z_1\cup Z_2$.}\\

\section{Automorphisms} 

\subsection{Finiteness of the Automorphism group}

Let $\pi\colon \tilde{X}\rightarrow X$ be the minimal resolution of singularities of a $4$-nodal tetrahedral quartic $X\subset {\mathbb P}^3$. Then  projections from each of the nodes give rise to $4$ involutions on $\tilde{X}$. We will use this observation to show that $Aut(\tilde{X})$ is an infinite group.\\

Let us denote by ${\phi}_i \colon \tilde{X}\rightarrow \tilde{X}$ the involution coming from the projection from node $E_i$ of $X$ and by ${\phi}_{i}^{*}\colon NS(\tilde{X})\rightarrow NS(\tilde{X})$ the corresponding automorphism of Neron-Severi group.\\

First, let us describe how ${\phi}_i$ acts on the elements of $M\subset NS(\tilde{X})$. Let us consider the case $i=4$.\\

{\bf Lemma 3.1.} {\it ${\phi}_{4}^{*}(L_{12})=R_3$, ${\phi}_{4}^{*}(L_{23})=R_1$, ${\phi}_{4}^{*}(L_{13})=R_2$, ${\phi}_{4}^{*}(L_{14})=E_1$, ${\phi}_{4}^{*}(L_{24})=E_2$, ${\phi}_{4}^{*}(L_{34})=E_3$, ${\phi}_{4}^{*}(R_3)=L_{12}$, ${\phi}_{4}^{*}(R_1)=L_{23}$, ${\phi}_{4}^{*}(R_2)=L_{13}$, ${\phi}_{4}^{*}(E_2)=L_{24}$, ${\phi}_{4}^{*}(E_3)=L_{34}$, ${\phi}_{4}^{*}(E_1)=L_{14}$, ${\phi}_{4}^{*}(R_4)=H-E_4-R_4$, ${\phi}_{4}^{*}(H)=R_1+R_2+2R_3-R_4+E_1+E_2+L_{12}+L_{34}+2(L_{14}+L_{24})$, ${\phi}_{4}^{*}(E_4)=2R_1-E_1+E_2+E_3-E_4+2L_{23}+L_{24}+L_{34}-L_{14}$. In particular, ${\phi}_{4}^{*}(M)\subset M$.}\\

{\it Proof:} Most of the equalities follow by continuity. Let us compute ${\phi}_{4}^{*}(R_4)$ and ${\phi}_{4}^{*}(E_4)$.\\

Let $C$ be the plane nodal rational curve on $X$ obtained by intersecting $X$ with the hyperplane in ${\mathbb P}^3$ passing through $E_4$ and $R_4$. Then its strict transform in $\tilde{X}$ is a smooth irreducible rational curve, whose class $H-E_4-R_4$ is exactly ${\phi}_{4}^{*}(R_4)$, because involution ${\phi}_4$ interchanges $C$ and $R_4$.\\

Since $E_4=H-E_1-E_2-R_3-L_{14}-L_{24}-L_{12}=E_3+R_4-R_3+L_{13}+L_{23}-L_{14}-L_{24}$, we have that ${\phi}_{4}^{*}(E_4)=$ $L_{34}+H-E_4-R_4-L_{12}+R_2+R_1-E_1-E_2=$ $2R_1-E_1+E_2+E_3-E_4+2L_{23}+L_{24}+L_{34}-L_{14}$. {\it QED}\\

{\bf Proposition 3.2.} {\it ${{\phi}_{4}^{*}}\circ {{\phi}_{3}^{*}}\in Aut(NS(\tilde{X})_{\mathbb Q})$ has infinite order.}\\

{\it Proof:} We will compute the matrix of the restriction of ${{\phi}_{4}^{*}}\circ {{\phi}_{3}^{*}}$ onto $M_{\mathbb Q}$ and check that it has infinite order.\\

The matrix of ${{\phi}_{4}^{*}}{\mid}_M$ in the basis $l_{12}, l_{13}, l_{14}, l_{23}, l_{24}, l_{34}, e_1, e_2, e_3, e_4, r_1$ is:
$$
\alpha=\left[
   \begin{array}{ccccccccccc}
-1 & 0 & 0 & 0 & 0 & 0 & 0 & 0 & 0 & 0 & 0\\
0 & -1 & 0 & 0 & 0 & 0 & 0 & 0 & 0 & 0 & 0\\
-1 & -1 & 0 & 0 & 0 & 0 & 1 & 0 & 0 & -1 & 0\\
1 & 1 & 0 & 0 & 0 & 0 & 0 & 0 & 0 & 2 & 1\\
0 & 1 & 0 & 0 & 0 & 0 & 0 & 1 & 0 & 1 & 0\\
1 & 0 & 0 & 0 & 0 & 0 & 0 & 0 & 1 & 1 & 0\\
-1 & -1 & 1 & 0 & 0 & 0 & 0 & 0 & 0 & -1 & 0\\
0 & 1 & 0 & 0 & 1 & 0 & 0 & 0 & 0 & 1 & 0\\
1 & 0 & 0 & 0 & 0 & 1 & 0 & 0 & 0 & 1 & 0\\
0 & 0 & 0 & 0 & 0 & 0 & 0 & 0 & 0 & -1 & 0\\
1 & 1 & 0 & 1 & 0 & 0 & 0 & 0 & 0 & 2 & 0
   \end{array}
 \right]
$$

(we recall that $R_2=H-E_1-E_3-E_4-L_{13}-L_{14}-L_{34}$ $=R_1+E_2-E_1+L_{23}+L_{24}-L_{13}-L_{14}$, $R_3=H-E_1-E_2-E_4-L_{12}-L_{14}-L_{24}$ $=R_1+E_3-E_1+L_{23}+L_{34}-L_{14}-L_{12}$ and $R_4=H-E_1-E_2-E_3-L_{12}-L_{13}-L_{23}$ $=R_1+E_4-E_1+L_{24}+L_{34}-L_{12}-L_{13}$)\\

Similarly, one computes the matrix of ${{\phi}_{3}^{*}}{\mid}_M$ in the same basis: 
$$
\beta=\left[
   \begin{array}{ccccccccccc}
-1 & 0 & 0 & 0 & 0 & 0 & 0 & 0 & 0 & 0 & 0\\
-1 & 0 & -1 & 0 & 0 & 0 & 1 & 0 & -1 & 0 & 0\\
0 & 0 & -1 & 0 & 0 & 0 & 0 & 0 & 0 & 0 & 0\\
0 & 0 & 1 & 0 & 0 & 0 & 0 & 1 & 1 & 0 & 0\\
1 & 0 & 1 & 0 & 0 & 0 & 0 & 0 & 2 & 0 & 1\\
1 & 0 & 0 & 0 & 0 & 0 & 0 & 0 & 1 & 1 & 0\\
-1 & 1 & -1 & 0 & 0 & 0 & 0 & 0 & -1 & 0 & 0\\
0 & 0 & 1 & 1 & 0 & 0 & 0 & 0 & 1 & 0 & 0\\
0 & 0 & 0 & 0 & 0 & 0 & 0 & 0 & -1 & 0 & 0\\
1 & 0 & 0 & 0 & 0 & 1 & 0 & 0 & 1 & 0 & 0\\
1 & 0 & 1 & 0 & 1 & 0 & 0 & 0 & 2 & 0 & 0
   \end{array}
 \right]
$$

Then 
$$
{\alpha}\circ{\beta}=\left[
   \begin{array}{ccccccccccc}
1 & 0 & 0 & 0 & 0 & 0 & 0 & 0 & 0 & 0 & 0\\
1 & 0 & 1 & 0 & 0 & 0 & -1 & 0 & 1 & 0 & 0\\
0 & 1 & 0 & 0 & 0 & -1 & -1 & 0 & -1 & 0 & 0\\
1 & 0 & 0 & 0 & 1 & 2 & 1 & 0 & 3 & 0 & 0\\
0 & 0 & 0 & 1 & 0 & 1 & 1 & 0 & 1 & 0 & 0\\
0 & 0 & 0 & 0 & 0 & 1 & 0 & 0 & 0 & 0 & 0\\
1 & 0 & 0 & 0 & 0 & -1 & -1 & 0 & 0 & 0 & 0\\
1 & 0 & 0 & 0 & 0 & 1 & 1 & 0 & 2 & 0 & 1\\
1 & 0 & 0 & 0 & 0 & 1 & 0 & 0 & 2 & 1 & 0\\
-1 & 0 & 0 & 0 & 0 & -1 & 0 & 0 & -1 & 0 & 0\\
0 & 0 & 0 & 0 & 0 & 2 & 1 & 1 & 2 & 0 & 0
   \end{array}
 \right]
$$

One checks by a direct computation that the matrix ${\alpha}\circ{\beta}$ is of infinite order. {\it QED}\\

{\bf Corollary 3.3.} {\it $Aut(\tilde{X})=Bir(X)$ is infinite.}\\

{\it Proof:} ${{\phi}_3}\circ{{\phi}_4}\in Aut(\tilde{X})$ has infinite order by the Propositions 3.2 above. {\it QED}\\

There is another way to construct involutions on tetrahedral quartics. Let $\pi\colon \tilde{X}\rightarrow X$ be the minimal resolution of singularities of a general tetrahedral quartic. Let $H_1=L_{12}+L_{13}+L_{14}+L_{23}+L_{24}+L_{34}+R_1+R_2+R_3+R_4$. It is immediate that $H_1$ is nef and is not hyperelliptic (see \cite{SD}). Consider the morphism ${{\phi}_{\mid H_1 \mid}}\colon \tilde{X}\rightarrow {\mathbb P}^3$ given by the complete linear system $\mid H_1 \mid$. It contracts the $6$ curves $L_{ij}$ into the $6$ nodes and represents $\tilde{X}$ as the minimal resolution of singularities of a quartic in ${\mathbb P}^3$ with exactly $6$ nodes (the images of $L_{ij}$) and exactly $4$ lines (the images of $R_i$) such that all the nodes lie in the same plane in ${\mathbb P}^3$, each line contains $3$ nodes and every node lies on $2$ lines.

\begin{center}
  \includegraphics{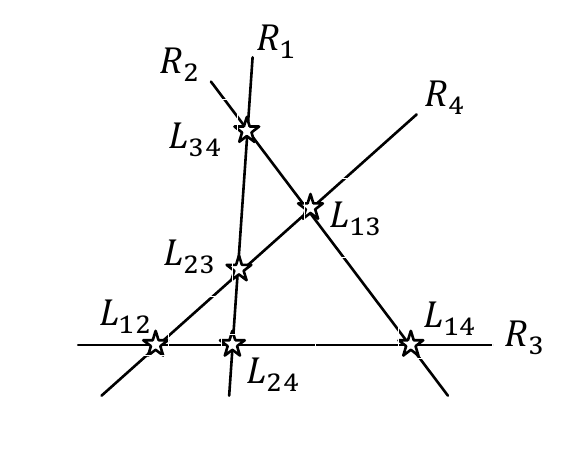}
\end{center} 

The images of the curves $E_i$ are twisted cubics, which intersect the plane in ${\mathbb P}^3$ at the corresponding $3$ nodes. Now, projection from each of the $6$ nodes gives an involution on $\tilde{X}$.\\

So, we see that to each of the $14$ curves $L_{ij}$, $E_i$, $R_j$ on $\tilde{X}$ there corresponds an involution on $\tilde{X}$.\\

{\bf Question.} Do these $14$ involutions generate $Aut(\tilde{X})$?

\subsection{The image of $Aut(\tilde{X})\rightarrow Aut(A_M,q_M)$}

Let us assume now that $X$ is a very general tetrahedral quartic, which in particular means that $NS(\tilde{X})\cong M$, where lattice $M$ is defined as follows: 
$$
M=\left[
   \begin{array}{ccccccccccc}
-2 & 0 & 0 & 0 & 0 & 0 & 1 & 1 & 0 & 0 & 0\\
0 & -2 & 0 & 0 & 0 & 0 & 1 & 0 & 1 & 0 & 0\\
0 & 0 & -2 & 0 & 0 & 0 & 1 & 0 & 0 & 1 & 0\\
0 & 0 & 0 & -2 & 0 & 0 & 0 & 1 & 1 & 0 & 1\\
0 & 0 & 0 & 0 & -2 & 0 & 0 & 1 & 0 & 1 & 1\\
0 & 0 & 0 & 0 & 0 & -2 & 0 & 0 & 1 & 1 & 1\\
1 & 1 & 1 & 0 & 0 & 0 & -2 & 0 & 0 & 0 & 0\\
1 & 0 & 0 & 1 & 1 & 0 & 0 & -2 & 0 & 0 & 0\\
0 & 1 & 0 & 1 & 0 & 1 & 0 & 0 & -2 & 0 & 0\\
0 & 0 & 1 & 0 & 1 & 1 & 0 & 0 & 0 & -2 & 0\\
0 & 0 & 0 & 1 & 1 & 1 & 0 & 0 & 0 & 0 & -2
   \end{array}
 \right]
$$

It has a discriminant-group $A_M=M^{*}/M$, a discriminant-form $q_M\colon A_M\rightarrow {\mathbb Q}/2{\mathbb Z}$ and a finite bilinear form $b_M\colon A_M\times  A_M\rightarrow {\mathbb Q}/{\mathbb Z}$. Let us denote by $l_{12}, l_{13}, l_{14}, l_{23}, l_{24}, l_{34}, e_1, e_2, e_3, e_4, r_1$ the basis in ${\mathbb Z}^{11}$, in which the intersection matrix of $M$ has the above form. Let us denote by ${\lambda}_{12}, {\lambda}_{13}, {\lambda}_{14}, {\lambda}_{23}, {\lambda}_{24}, {\lambda}_{34}, {\epsilon}_1, {\epsilon}_2, {\epsilon}_3, {\epsilon}_4, {\rho}$ the corresponding elements of the dual basis in $M^{*}=Hom_{\mathbb Z}(M,{\mathbb Z})$.\\

{\bf Lemma 3.4.} {\it $A_M={\mathbb Z}{\epsilon}_1 \oplus {\mathbb Z}{\lambda}_{23} \oplus {\mathbb Z}{\lambda}_{24}\cong$ ${\mathbb Z}/8{\mathbb Z} \oplus {\mathbb Z}/4{\mathbb Z} \oplus {\mathbb Z}/4{\mathbb Z}$.\\

$b_M({\epsilon}_1,{\epsilon}_1)=\frac{3}{8}+{\mathbb Z}$, $b_M({\lambda}_{23},{\lambda}_{23})=b_M({\lambda}_{24},{\lambda}_{24})=-\frac{1}{2}+{\mathbb Z}$, $b_M({\epsilon}_1,{\lambda}_{23}) = b_M({\epsilon}_1,{\lambda}_{24})=\frac{1}{2}+{\mathbb Z}$, $b_M({\lambda}_{23},{\lambda}_{24})=\frac{1}{4}+{\mathbb Z}$.\\

$q_M({\epsilon}_1)=\frac{3}{8}+2{\mathbb Z}$, $q_M({\lambda}_{23})=q_M({\lambda}_{24})=-\frac{1}{2}+2{\mathbb Z}$.}\\

{\it Proof:} Elementary. Note that in $A_M$ ${\lambda}_{12}=-2 {\epsilon}_1 -{\lambda}_{23}-{\lambda}_{24}$, ${\lambda}_{13}=2 {\epsilon}_1 +{\lambda}_{24}$, ${\lambda}_{14}=2 {\epsilon}_1 +{\lambda}_{23}$, ${\lambda}_{34}=4 {\epsilon}_1 -{\lambda}_{23}-{\lambda}_{24}$, ${\epsilon}_{2}=3 {\epsilon}_1 +2{\lambda}_{23}+2{\lambda}_{24}$, ${\epsilon}_{3}=3 {\epsilon}_1 +2{\lambda}_{24}$, ${\epsilon}_{4}=3 {\epsilon}_1 +2{\lambda}_{23}$, $\rho =2 {\epsilon}_1$.\\

Note also that in $M^{*}$:
$$
8\cdot{\epsilon}_1=-13e_1+e_2+e_3+e_4+6r_1-6l_{12}-6l_{13}-6l_{14}+4l_{23}+4l_{24}+4l_{34}
$$

$$
4\cdot{\lambda}_{23}=2e_1+2e_4+l_{12}+l_{13}+2l_{14}-2l_{23}+l_{24}+l_{34}
$$

$$
4\cdot{\lambda}_{24}=2e_1+2e_3+l_{12}+2l_{13}+l_{14}+l_{23}-2l_{24}+l_{34}
$$ {\it QED}\\

Now we will study the image of the composition of the natural group homomorphisms $Aut(\tilde{X})\rightarrow Aut(NS(\tilde{X}))=Aut(M)\rightarrow Aut(A_M)$. First of all, let us describe the images in $A_M$ of some natural subgroups of $Aut(M)$.\\

The symmetric group $S_4$ acts naturally on $M$ by permuting indices $\{ 1,2,3,4 \}$. It is immediate that the corresponding action of transpositions on $A_M$ is as follows:\\
$$
(12)({\epsilon}_1)={\epsilon}_1+2{\lambda}_{23}+2{\lambda}_{24}, \; (12)({\lambda}_{23})=4{\epsilon}_1+{\lambda}_{24}, \; (12)({\lambda}_{24})=4{\epsilon}_1+{\lambda}_{23},
$$

$$
(13)({\epsilon}_1)={\epsilon}_1+2{\lambda}_{24}, \; (13)({\lambda}_{23})=-{\lambda}_{23}-{\lambda}_{24}, \; (13)({\lambda}_{24})={\lambda}_{24},
$$

$$
(14)({\epsilon}_1)={\epsilon}_1+2{\lambda}_{23}, \; (14)({\lambda}_{23})={\lambda}_{23}, \; (14)({\lambda}_{24})=-{\lambda}_{23}-{\lambda}_{24},
$$

$$
(23)({\epsilon}_1)={\epsilon}_1, \; (23)({\lambda}_{23})={\lambda}_{23}, \; (23)({\lambda}_{24})=4{\epsilon}_1-{\lambda}_{23}-{\lambda}_{24},
$$

$$
(24)({\epsilon}_1)={\epsilon}_1, \; (24)({\lambda}_{23})=4{\epsilon}_1-{\lambda}_{23}-{\lambda}_{24}, \; (24)({\lambda}_{24})={\lambda}_{24},
$$

$$
(34)({\epsilon}_1)={\epsilon}_1, \; (34)({\lambda}_{23})={\lambda}_{24}, \; (34)({\lambda}_{24})={\lambda}_{23}.
$$

Interchanging nodes and residual lines of $X$ (i.e. mirror pairing) leads to an automorphism of $M$, and hence to an automorphism $\mu$ of $A_M$. A direct computation shows that:
$$
\mu({\epsilon}_1)={\epsilon}_1, \; \mu({\lambda}_{23})=-{\lambda}_{23}, \; \mu({\lambda}_{24})=-{\lambda}_{24}.
$$

Finally, the covering involution ${\phi}_4$ (coming from the projection of $X$ from one of its nodes) leads to the following involution $i$ on $A_M$:
$$
i({\epsilon}_1)=-{\epsilon}_1, \; i({\lambda}_{23})=-{\lambda}_{23}, \; i({\lambda}_{24})=-{\lambda}_{24}.
$$

{\bf Lemma 3.5.} {\it An automorphism of the finite group $A_M$ preserves the discriminant-form $q_M$, if and only if it is a composition of the automorphisms induced by $S_4$, $\mu$ and $i$.}\\

{\it Proof:} Direct verification. {\it QED}\\

{\bf Remark.} It follows from Nikulin's results \cite{Nikulin} that the group homomorphism $Aut(M)\rightarrow Aut(A_M,q_M)$ is surjective, since $rk(M)=11$ and $l(A_M)=l((A_M)_2)=3$ (in \cite{Nikulin} $l(G)$ denotes the minimal number of generators of a finite group $G$).\\ 

{\bf Corollary 3.6.} {\it The image of the composition $Aut(\tilde{X})\rightarrow Aut(NS(\tilde{X}))=Aut(M)\rightarrow Aut(A_M)$ consists of two elements: identity and $i$.}\\

{\it Proof:} Identity and $i$ lie in this image by construction (they both are induced by automorphisms of $\tilde{X}$).\\ 

Let $\alpha = {\sigma}\circ {\mu}$ or $\alpha = {\sigma}$ be an automorphism of $A_M$, where ${\sigma}\in S_4$, and let $\bar{\alpha}$ be the automorphism of $M$, which induces $\alpha$ on $A_M$. Suppose that there exists an automorphism $\phi$ of $\tilde{X}$, which induces $\alpha$ on $A_M$. This implies that there exists a Hodge isometry $\beta$ of the transcendental lattice $N$ of $\tilde{X}$ (i.e. the orthogonal complement of $M=NS(\tilde{X})$ in $H^2(\tilde{X},\mathbb Z)$), which induces the same automorphism $\alpha$ on $A_N\cong A_M$.\\

Since the image under $\bar{\alpha}$ of an ample divisor (say, $A=L_{12}+L_{13}+L_{14}+L_{23}+L_{24}+L_{34}+E_1+E_2+E_3+E_4+R_1+R_2+R_3+R_4$) is ample, it follows that we can glue $\bar{\alpha}$ and $\beta$ together into an effective Hodge isometry on $H^2(\tilde{X},\mathbb Z)$. Torelli theorem for $K3$ surfaces \cite{PS} implies that there exists an automorphism $\psi$ of $\tilde{X}$, which induces $\bar{\alpha}$ on $NS(\tilde{X})=M$. By definition of $\bar{\alpha}$ this automorphism permutes the curves $\{  L_{12}, L_{13}, L_{14}, L_{23}, L_{24}, L_{34}, E_1, E_2, E_3, E_4, R_1, R_2, R_3, R_4 \}$ among themselves. Corollary 0.7 implies that $\psi$ is the identity (since $X$ is general). Hence $\bar{\alpha}$ is the identity, and $\alpha$ is the identity. {\it QED}\\

{\bf Corollary 3.7.} {\it Let $X\subset {\mathbb P}^3$ be a very general tetrahedral quartic with $NS(\tilde{X})\cong M$, where $\pi\colon \tilde{X}\rightarrow X$ is the minimal resolution of singularities, and ${\pi}^{v}\colon \tilde{X}\rightarrow X^{v}$ is the minimal resolution of singularities of the mirror dual $X^{v}$ of $X$.\\

If ${\pi}_1\colon \tilde{X}\rightarrow {\mathbb P}^3$ represents $\tilde{X}$ as the minimal resolution of singularities of a tetrahedral quartic $Y$ in ${\mathbb P}^3$, then (upto an automorphism of ${\mathbb P}^3$) either ${\pi}_1={\pi}$, or ${\pi}_1={\pi}^{v}$.}\\

{\bf Remark.} In other words, mirror duality is indeed a pairing (at least, for very general tetrahedral quartics): only two (very general) tetrahedral quartics share the same minimal resolution of singularities - $X$ and its mirror $X^{v}$.\\

{\it Proof:} ${\pi}_1$ induces an identification $\alpha \colon NS(\tilde{X})\cong M$, and hence an automorphism of the lattice $M$, which sends $A\in M$ to an ample divisor class. By reordering vertices of the tetrahedron (corresponding to the representation ${\pi}_1\colon \tilde{X}\rightarrow {\mathbb P}^3$) and by applying the mirror construction to $Y$ we can ensure that the automorphism of $A_M$ induced by $\alpha$ is either identity or $i$. In either case, $\alpha$ can be glued with a suitable Hodge isometry on the transcendental lattice $N$ of $\tilde{X}$ to obtain an effective Hodge isometry on $H^2(\tilde{X},\mathbb Z)$, which by Torelli theorem comes from an automorphism $\gamma$ of $\tilde{X}$ and which by construction induces $\alpha$ on $NS(\tilde{X})\cong M$. Since both $\pi$ and ${\pi}_1$ are induced by complete linear systems (which are identified by $\alpha$), this implies that ${{\pi}_1}\circ \gamma=\pi$. This also implies that there exists an automorphism ${\gamma}_1$ of ${\mathbb P}^3$ such that ${{\gamma}_1}\circ {\pi}={\pi}_1$. {\it QED}\\

One can also use the study of $A_M$ shown above in order to give a certain description of the automorphism group of the minimal resolution of singularities $\tilde{X}$ of a very general tetrahedral quartic $X\subset {\mathbb P}^3$ (such that $NS(\tilde{X})\cong M$). Let $N$ be the transcendental lattice of $\tilde{X}$ (i.e. the orthogonal complement of $M$ in $H^2(\tilde{X},\mathbb Z)$). Recall that by the construction of the period domain of $K3$ surfaces, the period of $\tilde{X}$ lies in ${\mathbb P}(N_{\mathbb C})$.\\

{\bf Corollary 3.8.} {\it Let $\pi\colon \tilde{X}\rightarrow X$ is the minimal resolution of singularities of a very general tetrahedral quartic $X$ (such that $NS(\tilde{X})\cong M$). Let $N$ be the transcendental lattice of $K3$ surface $\tilde{X}$ (i.e. the orthogonal complement of $NS(\tilde{X})\cong M$ in $H^2(\tilde{X},\mathbb Z)$) together with the Hodge structure inherited from the Hodge structure of $H^2(\tilde{X},\mathbb Z)$. Then any Hodge isometry of $N$ is either identity or scalar multiplication by $-1$. It is $-1$ if and only if it induces automorphism $i$ on $A_N\cong A_M$.}\\

{\bf Remark.} In other words, the natural group homomorphism $Aut(N)\rightarrow Aut(A_N)$ induces the isomorphism between the subgroup of $Aut(N)$ consisting of those automorphisms of $N$, which preserve (upto scalar multiplication) the period of $\tilde{X}$ and the subgroup of $Aut(A_N)\cong Aut(A_M)$ consisting of two elements $\{ 1,i \}\cong {\mathbb Z}/2{\mathbb Z}$.\\

{\it Proof:} Since both identity and $i$ are induced (as automorphisms of $A_N\cong A_M$) by automorphisms of $\tilde{X}$, they are also induced by the Hodge isometries on $N$. Let $\beta$ be a Hodge isometry on $N$ and $\alpha$ be the automorphism of $A_N\cong A_M$ which is induced by $\beta$. If $\alpha$ were neither identity, nor $i$, we would be able to take a nontrivial (i.e. not equal to identity) automorphism $\bar{\alpha}$ on $M$ of the form $\sigma$ or ${\sigma}\circ {\mu}$ (where $\sigma \in S_4$), which induces $\alpha$ on $A_M$, and which sends $A\in M$ into an ample divisor. Then we would be able to glue $\bar{\alpha}$ and $\beta$ into an effective Hodge isometry on $H^2(\tilde{X},\mathbb Z)$. So, by Torelli theorem we would obtain an automorphism of $\tilde{X}$, which induces on $A_N$ an automorphism different from $i$ and identity. This contradicts to Corollary 3.6. Hence the image of the group of Hodge isometries of $N$ in the automorphism group of $A_N\cong A_M$ (under the restriction homomorphism) consists of two elements - identity and $i$.\\

Now let us show that if a Hodge isometry $\beta$ on $N$ induces the identity automorphism on $A_N\cong A_M$, then $\beta$ is also identity. If it were not, than by Torelli theorem the effective Hodge isometry of $H^2(\tilde{X},\mathbb Z)$, which induces identity on $M$ and $\beta$ on $N$ is induced by an automorphism of $\tilde{X}$. This automorphism should be identity by Corollary 0.7. Hence $\beta$ is identity. {\it QED}\\

{\bf Corollary 3.9.} {\it Let $G_M$ be the subgroup of $Aut(M)$, consisting of those isometries of $M$, which induce identity on $A_M$ and which send $A\in M\cong NS(\tilde{X})$ into an ample divisor class. Then there is a right split short exact sequence of groups:
$$
0\rightarrow G_M \rightarrow Aut(\tilde{X}) \rightarrow {\mathbb Z}/2{\mathbb Z} \rightarrow 0
$$
}\\

{\it Proof:} The existence of the epimorphism $Aut(\tilde{X})\rightarrow {\mathbb Z}/2{\mathbb Z} \cong \{ 1,i \}\hookrightarrow Aut(A_M)$ was shown above. The right splitting is given by the projection from any node of $X$ (we saw above that such a projection induces automorphism $i$ on $A_M$). {\it QED}\\

{\bf Remark.} The condition that an element of $G_M$ preserves ampleness of $A\in M$ can be reformulated as follows. If $\alpha \in ker(Aut(M)\rightarrow Aut(A_M))$, then $\alpha$ lies in $G_M$ if and only if for any $c\in M$ such that $c^2=-2$ and $c\cdot A>0$ we have that ${\alpha}(A)\cdot c > 0$. Note that, it does not depend on $X$.\\

It is also easy to see that $G_M$ as a subgroup of $Aut(\tilde{X})$ is exactly the subgroup of symplectic automorphisms of $\tilde{X}$. (An automorphism of a $K3$ surface $Y$ is called symplectic, if the Hodge isometry induced by it on the second cohomology group $H^2(Y,\mathbb Z)$ preserves the period $\omega={\omega}(Y)\in H^2(Y,\mathbb Z)\otimes \mathbb C$.)\\ 

{\bf Corollary 3.10.} {\it If $X\subset {\mathbb P}^3$ is a very general tetrahedral quartic such that $NS(\tilde{X})\cong M$ (where $\pi\colon \tilde{X}\rightarrow X$ is the minimal resolution of singularities), then the automorphism group of $\tilde{X}$ does not depend on $X$ and is isomorphic to a semidirect product:
$$
Aut(\tilde{X})\cong G_M  \rtimes ({\mathbb Z}/2{\mathbb Z}).
$$
Group $G_M$ depends only on the lattice $M$ and consists of those $\alpha \in ker(Aut(M)\rightarrow Aut(A_M))$ such that ${\alpha}(A)\cdot c>0$ for any $c \in M$ satisfying $A\cdot c>0$ and $c^2=-2$.\\
The generator of ${\mathbb Z}/2{\mathbb Z}$ acts on $G_M$ by sending $\alpha \in G_M$ to ${{\phi}_4^{*}}\circ{\alpha}\circ {{\phi}_4^{*}}$, where ${{\phi}_4^{*}}$ is the isometry of $M$ induced by the covering involution corresponding to the projection of the tetrahedral quartic $X$ from one of its (chosen) nodes.}\\

\section{Double del Pezzo surfaces} 

We have already observed in the Introduction that if $\pi\colon \tilde{X}\rightarrow X$ is the minimal resolution of singularities of a general (in particular, $4$-nodal) tetrahedral quartic $X\subset {\mathbb P}^3$, then projection from a node of $X$ represents $\tilde{X}$ as the minimal resolution of singularities of the double cover of ${\mathbb P}^2$ ramified over a sextic curve $C(X)\subset {\mathbb P}^2$.\\

{\bf Lemma 4.1.} {\it For a general tetrahedral quartic $X$, $C(X)$ is a plane irreducible sextic with $3$ cusps at the vertices of a (nondegenerate) triangle (and no other singularities), which has a tritangent line $L$ (which touches $C(X)$ at $3$ distinct smooth points) and a (smooth) tritangent conic $Q$ which passes through the cusps (and touches $C(X)$ at $3$ distinct smooth points).\\

Moreover, each edge of the triangle touches $C(X)$ at one smooth point, and there is a (smooth) plane cubic passing through the cusps of $C(X)$, the $9$ points of tangency of $C(X)$ with $Q$, $L$ and the edges of the triangle and having at the cusps of $C(X)$ the same tangent lines as $C(X)$.}\\

{\it Proof:} Elementary. \\If the equation of $X\subset {\mathbb P}^3$ is $F(X_0,X_1,X_2,X_3)=0$, where $F(X_0,X_1,X_2,X_3)=$ $A(X_0,X_1,X_2)\cdot X_0X_1X_2+$ $B(X_0,X_1,X_3)\cdot X_0X_1X_3+$ $C(X_0,X_2,X_3)\cdot X_0X_2X_3+D(X_1,X_2,X_3)\cdot X_1X_2X_3+ X_0X_1X_2X_3$, $A(X_0,X_1,X_2)=a_0X_0+a_1X_1+a_2X_2$, $B(X_0,X_1,X_3)=b_0X_0+b_1X_1+b_3X_3$, $C(X_0,X_2,X_3)=c_0X_0+c_2X_2+c_3X_3$, $D(X_1,X_2,X_3)=d_1X_1+d_2X_2+d_3X_3$, then the equation of $C(X)\subset {\mathbb P}^2$ (corresponding to the projection from the node $(0:0:0:1)$ onto the opposte face of the tetrahedron) is $G(X_0,X_1,X_2)=0$, where $G(X_0,X_1,X_2)=[X_0X_1(b_0X_0+b_1X_1)+X_0X_2(c_0X_0+c_2X_2)+X_1X_2(d_1X_1+d_2X_2)+X_0X_1X_2]^2-4\cdot A(X_0,X_1,X_2)\cdot X_0X_1X_2\cdot (b_3 \cdot X_0X_1+c_3\cdot X_0X_2+d_3\cdot X_1X_2)$.\\

Line $L$ has equation $A(X_0,X_1,X_2)=0$ (and coincides with the residual line opposite to the node, from which $X$ was projected). The tritangent conic $Q$ has equation $b_3 \cdot X_0X_1+c_3\cdot X_0X_2+d_3\cdot X_1X_2=0$ (and coincides with the image of the node from which we project). One sees that each of the edges of the triangle touches $C(X)$ at a smooth point. The cubic touching $C(X)$ at the cusps and passing through its $9$ points of tangency with the edges of the triangle and the curves $Q$, $L$ is given by the equation $X_0X_1(b_0X_0+b_1X_1)+X_0X_2(c_0X_0+c_2X_2)+X_1X_2(d_1X_1+d_2X_2)+X_0X_1X_2=0$. {\it QED}\\

{\bf Remark.} Note that the tritangent conic comes together with $6$ distinguished points on it (the points of intersection of $Q$ with $C(X)$), and so naturally leads to a Kummer surface.\\

Blow-up of ${\mathbb P}^2$ at the cusps of $C(X)$ resolves singularities of $C(X)$ and represents $\tilde{X}$ as a double cover of the del Pezzo surface of degree $6$ ramified over the strict transform (i.e. the normalization) $\tilde{C}(X)$ of the sextic $C(X)$. We will see that vice-versa any such double cover of the del Pezzo surface of degree $6$ is (the minimal resolution of singularities of) a tetrahedral quartic in ${\mathbb P}^3$.\\

Let $C(X)\subset {\mathbb P}^2$ be a sextic as above (coming from a general tetrahedral quartic $X\subset {\mathbb P}^3$), $q\colon Z\rightarrow {\mathbb P}^2$ be the blow-up of ${\mathbb P}^2$ at the $3$ cusps of $C(X)$ and ${\phi}\colon Z\hookrightarrow {\mathbb P}^6$ be the closed immersion of $Z$ as a del Pezzo surface of degree $6$.\\

\begin{center}
  \includegraphics{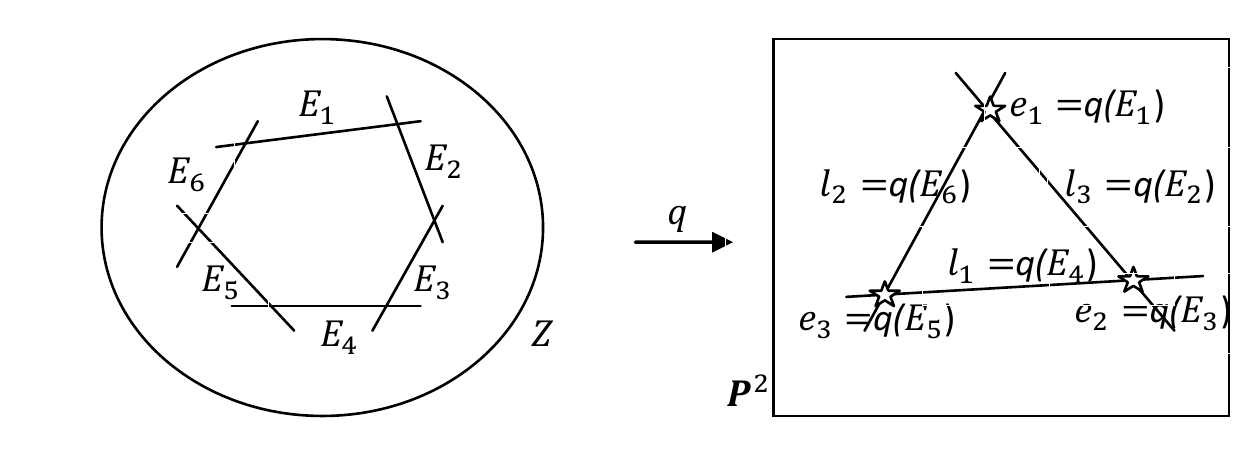}
\end{center}

Let $E_1, E_2, E_3, E_4, E_5, E_6$ be the $6$ lines on $Z$ such that $q(E_1)=e_1$, $q(E_3)=e_2$ and $q(E_5)=e_3$ are the cusps of $C(X)$ and $q(E_2)=l_3$, $q(E_4)=l_1$, $q(E_6)=l_2$ are the lines connecting them (the sides of the triangle). Let $\tilde{C}(X)$ be the strict transform of $C(X)$ in $Z$, i.e. $q {\mid}_{\tilde{C}(X)}\colon \tilde{C}(X)\rightarrow C(X)$ is the normalization of $C(X)$. We will use the same notation $L$ and $Q$ both for the tritangent line and conic in ${\mathbb P}^2$ and for their strict transforms (also tritangent to $\tilde{C}(X)$) in $Z$.\\

\begin{center}
  \includegraphics{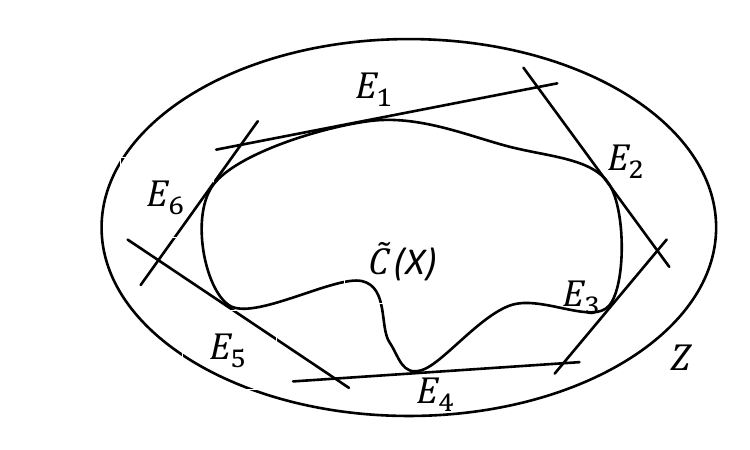}
\end{center}

Note that each of the $6$ lines $E_1, E_2, E_3, E_4, E_5, E_6$ touches $\tilde{C}(X)$ at one point. The $12$ points of tangency of $\tilde{C}(X)$ with $E_1, E_2, E_3, E_4, E_5, E_6$, $Q$, $L$ lie on a hyperplane in ${\mathbb P}^6$ (which intersects $\tilde{C}(X)$ exactly at these $12$ points).\\

The closed immersion ${\phi}\colon Z\rightarrow {\mathbb P}^6$ is given by the complete linear system (anticanonical) $\mid 3q^{*}h-e_1-e_2-e_3 \mid$ of cubics in ${\mathbb P}^2$ passing through the vertices of the triangle (where $h$ is the class of a line in ${\mathbb P}^2$). The class of $\tilde{C}(X)$ on $Z$ is $6q^{*}h-2e_1-2e_2-2e_3$. The class of the tritangent line $L$ on $Z$ is $q^{*}h$, and the class of the tritangent conic $Q$ on $Z$ is $2q^{*}h-e_1-e_2-e_3$.\\

Hence the composition $\tilde{C}(X)\hookrightarrow Z \xrightarrow{\phi} {\mathbb P}^6$ is the canonical embedding of the genus $7$ smooth curve $\tilde{C}(X)$. Note that the image of this composition is not degenerate, because $\mid (3q^{*}h-e_1-e_2-e_3)-(6q^{*}h-2e_1-2e_2-2e_3)  \mid=\emptyset$. In particular, $\tilde{C}(X)$ is not hyperelliptic.\\

We also see that $\tilde{C}(X)$ is cut out on $Z\subset {\mathbb P}^6$ by a quadric hypersurface $\Omega$ (since the restriction homomorphism $H^0({\mathbb P}^6, O(2))\rightarrow H^0(Z,-2K_Z)$ is surjective), which is tangent to the lines $E_1, E_2, E_3, E_4, E_5, E_6$ and to each of the twisted cubics $L$ and $Q$ (the images of the tritangent line and conic in ${\mathbb P}^6$) at $3$ distinct points.\\

Note that $\tilde{C}(X)$ has $3$ $g^{1}_{4}$ (coming from projections of $C(X)$ from each of the $3$ cusps) and $2$ $g^{2}_{6}$ (one of them corresponds to the morphism $q {\mid}_{\tilde{C}(X)}\colon \tilde{C}(X)\rightarrow {\mathbb P}^2$ and the other one is cut out by conics passing through the cusps of $C(X)$).\\ 

Note that these two $g^{2}_{6}$ are 'dual' to each other. Namely, if we contract $E_2, E_4, E_6$ by $q'\colon Z\rightarrow {\mathbb P}^2$, then the $g^{2}_{6}$ coming from $q'$ will be exactly the $g^{2}_{6}$ coming from conics via $q$ and the $g^{2}_{6}$ coming from $q$ will be exactly the $g^{2}_{6}$ coming from conics via $q'$. This duality reflects the mirror pairing of general tetrahedral quartics mentioned earlier (if $q$ corresponds to the projection of $X$ from node $E_i$, then $q'$ corresponds to the projection of $X^{v}$ from its node $R_i$).\\

\begin{center}
  \includegraphics{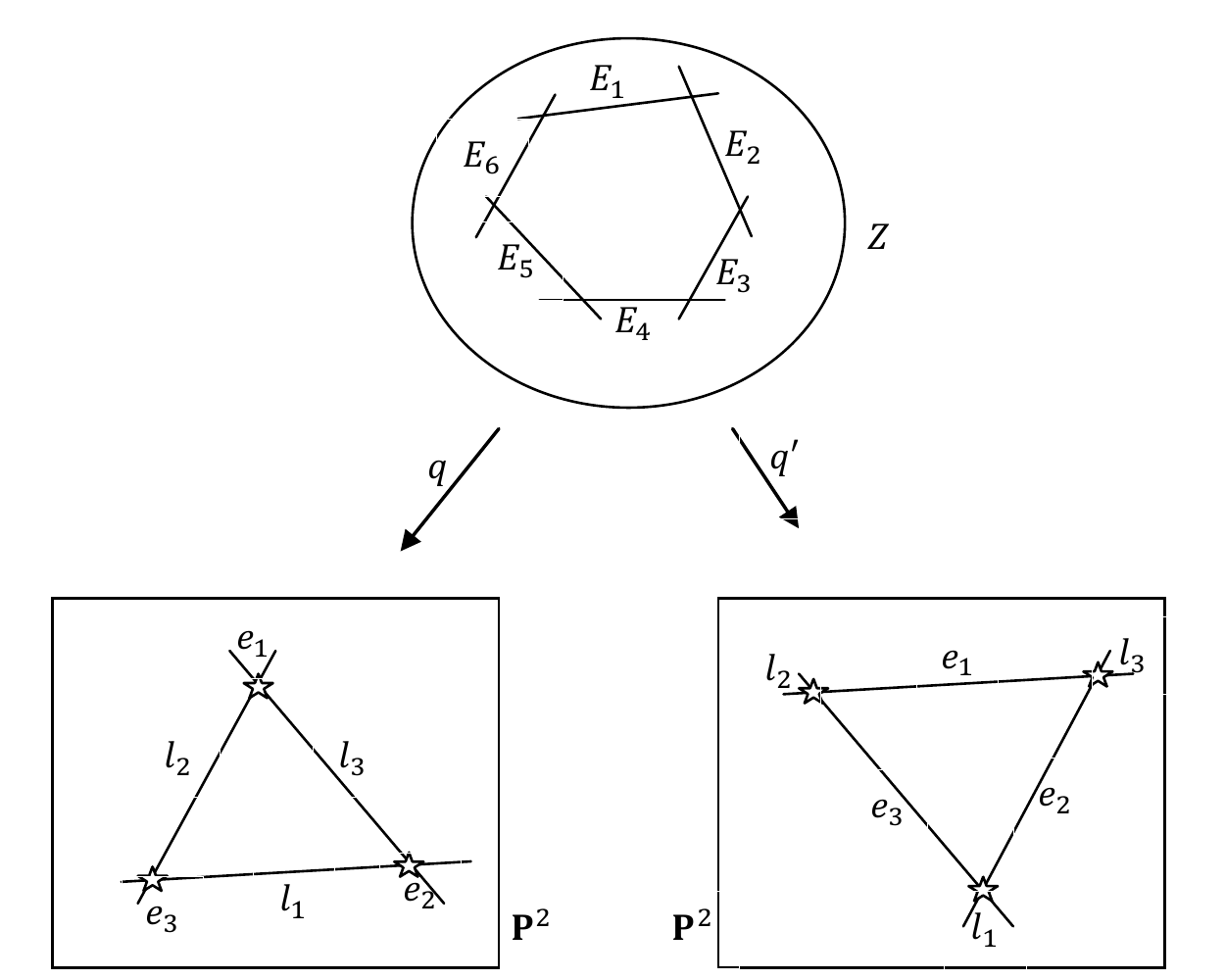}
\end{center}

We will call these $3$ $g^{1}_{4}$ and $2$ $g^{2}_{6}$ 'apparent'.\\

{\bf Lemma 4.2.} {\it If $X$ is a general tetrahedral quartic in ${\mathbb P}^3$, then the genus $7$ curve $\tilde{C}(X)$ has no $g^{1}_{2}$, has no $g^{1}_{3}$, has exactly $3$ $g^{1}_{4}$ (the 'apparent' ones), has no $g^{2}_{5}$ and has exactly $2$ $g^{2}_{6}$ (the 'apparent' ones).}\\  

{\it Proof:} (1) We have already observed that $\tilde{C}(X)$ has no $g^{1}_{2}$.\\

(2) Let $\delta$ be a $g^{1}_{3}$ on $\tilde{C}(X)$. Since $\tilde{C}(X)$ has no $g^{1}_{2}$, $\delta$ is base point free (and is a complete linear system by Clifford's theorem). Let $D=p_1+p_2+p_3$ be an element of $\delta$, where $p_1,p_2,p_3$ are $3$ distinct points on $\tilde{C}(X)$ not lying on the lines $E_1, E_2, E_3, E_4, E_5, E_6$.\\

By the geometric version of the Riemann-Roch theorem, $p_1$, $p_2$ and $p_3$ should lie on a line in ${\mathbb P}^6$. This means that $7$ distinct points $e_1, e_2, e_3, q(p_1), q(p_2), q(p_3), r$ (where $r$ is a general point on $C(X)$) in ${\mathbb P}^2$ fail to impose independent conditions on plane cubics. Hence $5$ of them should lie on a line (see Problem A-13, Chapter V, \cite{ACGH}). This line can not contain point $r$, because $r$ is general and $C(X)$ is irreducible. It can not contain two cusps of $C(X)$ either, because an edge of the triangle in ${\mathbb P}^2$ determined by the cusps of $C(X)$ intersects $C(X)$ only at one point ouside of the cusps. Hence this is impossible, and so $\tilde{C}(X)$ has no $g^{1}_{3}$.\\

(3) Let $\delta$ be a $g^{1}_{4}$ on $\tilde{C}(X)$. Since $\tilde{C}(X)$ has no $g^{1}_{3}$, $\delta$ is base point free (and is a complete linear system by Clifford's theorem). Let $D=p_1+p_2+p_3+p_4$ be an element of $\delta$, where $p_1,p_2,p_3,p_4$ are $4$ distinct points on $\tilde{C}(X)$ not lying on the lines $E_1, E_2, E_3, E_4, E_5, E_6$.\\

By the geometric version of the Riemann-Roch theorem, $p_1$, $p_2$, $p_3$ and $p_4$ span a $2$-plane in ${\mathbb P}^6$. This means that $7$ distinct points $e_1, e_2, e_3, q(p_1), q(p_2), q(p_3), q(p_4)$ in ${\mathbb P}^2$ fail to impose independent conditions on plane cubics. Hence $5$ of them should lie on a line (see Problem A-13, Chapter V, \cite{ACGH}). This line can not contain two cusps of $C(X)$, because such a line intersects $C(X)$ at one smooth point only. Hence $q(p_1), q(p_2), q(p_3), q(p_4)$ lie on a line in ${\mathbb P}^2$ passing through one of the cusps of $C(X)$. This means that $\delta$ is one of the $3$ 'apparent' $g^{1}_{4}$ on $\tilde{C}(X)$ .\\

(4) Let $\delta$ be a $g^{2}_{5}$ on $\tilde{C}(X)$. Since $\tilde{C}(X)$ has no $g^{1}_{3}$, $\delta$ is base point free (and is a complete linear system by Clifford's theorem). Let $D=p_1+p_2+p_3+p_4+p_5$ be an element of $\delta$, where $p_1,p_2,p_3,p_4,p_5$ are $5$ distinct points on $\tilde{C}(X)$ not lying on the lines $E_1, E_2, E_3, E_4, E_5, E_6$.\\

By the geometric version of the Riemann-Roch theorem, $p_1$, $p_2$, $p_3$, $p_4$, $p_5$ span a $2$-plane in ${\mathbb P}^6$. Suppose this $2$-plane is determined by $p_1$, $p_2$ and $p_3$. Then this implies that $7$ distinct points $e_1, e_2, e_3, q(p_1), q(p_2), q(p_3), q(p_4)$ in ${\mathbb P}^2$ fail to impose independent conditions on plane cubics. Hence (by Problem A-13, Chapter V, \cite{ACGH}) $q(p_1), q(p_2), q(p_3), q(p_4)$ lie on a line in ${\mathbb P}^2$ passing through a cusp of $C(X)$. For the same reason, points $q(p_1), q(p_2), q(p_3), q(p_5)$ lie on (the same) line in ${\mathbb P}^2$. Hence all the $5$ distinct points on $C(X)$ $q(p_1), q(p_2), q(p_3), q(p_4), q(p_5)$ should lie on a lines in ${\mathbb P}^2$ passing through a cusp of $C(X)$. This contradicts to Bezout's theorem. Hence $\tilde{C}(X)$ has no $g^{2}_{5}$.\\ 

(5) Let $\delta$ be a $g^{2}_{6}$ on $\tilde{C}(X)$. Since $\tilde{C}(X)$ has no $g^{2}_{5}$, $\delta$ is base point free (and is a complete linear system by Clifford's theorem). Let $D=p_1+p_2+p_3+p_4+p_5+p_6$ be an element of $\delta$, where $p_1,p_2,p_3,p_4,p_5,p_6$ are $6$ distinct points on $\tilde{C}(X)$ not lying on the lines $E_1, E_2, E_3, E_4, E_5, E_6$.\\

By the geometric version of the Riemann-Roch theorem, $p_1$, $p_2$, $p_3$, $p_4$, $p_5$, $p_6$ span a $3$-plane in ${\mathbb P}^6$. Suppose that this $3$-plane is determined by $p_1$, $p_2$, $p_3$ and $p_4$. Then this implies that $8$ distinct points $e_1, e_2, e_3, q(p_1), q(p_2), q(p_3), q(p_4), q(p_5)$ in ${\mathbb P}^2$ fail to impose independent conditions on plane cubics. Hence (by Problem A-14, Chapter V, \cite{ACGH}) either $4$ of the points $q(p_1), q(p_2), q(p_3), q(p_4), q(p_5)$ lie on a line in ${\mathbb P}^2$ passing through a cusp of $C(X)$, or all $5$ points $q(p_1), q(p_2), q(p_3), q(p_4), q(p_5)$ lie on a smooth conic passing through all $3$ cusps of $C(X)$, or $q(p_1), q(p_2), q(p_3), q(p_4), q(p_5)$ lie on a line in ${\mathbb P}^2$. The first option leads to a contradiction with Besout's theorem (since $q(p_6)$ would have to lie on the same line). Hence $q(p_1), q(p_2), q(p_3), q(p_4), q(p_5)$ lie either on a line in ${\mathbb P}^2$ or on a smooth conic in ${\mathbb P}^2$ passing through the $3$ cusps of $C(X)$. The same dilemma (with the same line and conic) holds for $q(p_1), q(p_2), q(p_3), q(p_4), q(p_6)$.\\

Hence all the $6$ points $q(p_1), q(p_2), q(p_3), q(p_4), q(p_5), q(p_6)$ in ${\mathbb P}^2$ lie either on a line or on a conic passing through all $3$ cusps of $C(X)$. Hence $\delta$ is one of the $2$ 'apparent' $g^{2}_{6}$ on $\tilde{C}(X)$. {\it QED}\\

Choose a triangle with vertices $e_1, e_2, e_3$ and edges $l_1, l_2, l_3$ in ${\mathbb P}^2$, a line $L$ and a smooth conic $Q$ such that $Q$ passes through the vertices of the triangle and $L$ does not pass through $e_1, e_2, e_3$.\\

\begin{center}
  \includegraphics{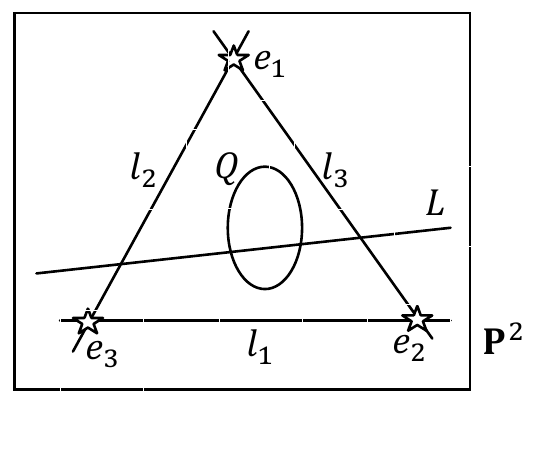}
\end{center}

Consider the degree $6$ del Pezzo surface $Z\subset {\mathbb P}^6$ obtained by blowing up vertices of the triangle.\\

Let $\Omega$ be a quadric in ${\mathbb P}^6$, which is neither a pair of planes nor a double plane, which touches all $6$ lines $E_1, E_2, E_3, E_4, E_5, E_6$ outside of their points of intersection and each of the $2$ twisted cubics $L$ and $Q$ (the images of the chosen line and conic in ${\mathbb P}^6$) at $3$ distinct points outside of $E_1, E_2, E_3, E_4, E_5, E_6$ and $Q\cap L$ in such a way that the $12$ points of tangency of $\Omega$ with $E_1, E_2, E_3, E_4, E_5, E_6, Q, L$ lie on a hyperplane in ${\mathbb P}^6$, and such that the curve of intersection $C={\Omega}\cap {Z}$ is smooth, irreducible and nondegenerate. (Alternatively, one may choose a cubic in ${\mathbb P}^2$, which touches neither of the curves $E_1, E_2, E_3, E_4, E_5, E_6, Q, L$ and require that $\Omega$ touch $E_1, E_2, E_3, E_4, E_5, E_6, Q, L$ at the points of their intersection with the cubic.)\\

{\bf Theorem 4.3.} {\it Let $r\colon \tilde{Z}\rightarrow Z$ be the double cover ramified over $C={\Omega}\cap {Z}$. Then $\tilde{Z}$ is the minimal resolution of singularities of a $4$-nodal tetrahedral quartic $X\subset {\mathbb P}^3$, $C=\tilde{C}(X)$ is the normalization of the branch sextic of the projection of $X$ from one of its nodes, and the double cover $r\colon \tilde{Z}=\tilde{X}\rightarrow Z$ comes from this projection.}\\  

{\it Proof:} This follows immediately from the comparison of equations of a tetrahedral quartic and of the branch sextic. We give here a more invariant argument.\\

It is immediate that $\tilde{Z}$ is a smooth $K3$ surface. Let $q\in {\Gamma}(Z, -2K_Z)$ and $c\in {\Gamma}(Z, -K_Z)$ be the sections corresponding to $C={\Omega}\cap {Z}$ and to the hyperplane in ${\mathbb P}^6$ containing the points of tangency of $\Omega$ with $E_1, E_2, E_3, E_4, E_5, E_6, Q, L$ respectively. Then there exists ${\epsilon}\in {\mathbb C}^{*}$ such that $q-({\epsilon}\cdot c)^2={{\lambda}_1}{{\lambda}_2}{{\lambda}_3}\cdot {\alpha}\cdot {\beta}$ in ${\Gamma}(Z, -2K_Z)={\Gamma}(Z,6q^{*}h-2e_1-2e_2-2e_3)\subset {\Gamma}({\mathbb P}^2, O_{{\mathbb P}^2}(6))$, where ${{\lambda}_1}, {{\lambda}_2}, {{\lambda}_3}, {\alpha}$ are linear forms on ${\mathbb P}^2$ corresponding to the lines $l_1, l_2, l_3, L$ respectively, and $\beta$ is a quadratic form on ${\mathbb P}^2$ corresponding to the conic $Q$.\\

Observe that $r^{*}(L)=L+L'$, $r^{*}(Q)=Q+Q'$, $r^{*}(E_i)=E_i+{E_i}'$ are pairs of smooth irreducible rational curves, which intersect transversely at (the preimages of) the points of intersection of $C$ with $L$, $Q$ and $E_i$ respectively. Let $\tilde{C}\subset \tilde{Z}$ be the reduced preimage of $C\subset Z$. We can choose notation in such a way that $L\cap {{E_i}'}=Q'\cap {{E_i}'}=L\cap Q=L'\cap Q'=L'\cap {{E_i}}=Q\cap {{E_i}}=\emptyset$ for all $i$ and ${{E_i}'}\cap {{E_j}'}=\emptyset$, if $i\neq j$.\\

Indeed, if we denote $T=L\cup Q \cup E_1 \cup E_2 \cup E_3 \cup E_4 \cup E_5 \cup E_6$, then the restriction of the double covering $r\colon \tilde{Z}\rightarrow Z$ onto $T\subset Z$ has a section, whose image in $\tilde{Z}=Spec({{\mathcal O}_Z}\oplus {K_Z})=Spec(\frac{{{\mathcal O}_Z}[t]}{t^2-q})$ may be taken to be $div(t-{\epsilon}c)\subset \tilde{Z}$, because $r^{-1}(T)=Spec(\frac{{{\mathcal O}_Z}[t]}{t^2-({\epsilon}c)^2})=Spec(\frac{{{\mathcal O}_Z}[t]}{t-{\epsilon}c})\cup Spec(\frac{{{\mathcal O}_Z}[t]}{t+{\epsilon}c})$.\\ 

Let us also note that $Q'+L+E_1+E_2+E_3+E_4+E_5+E_6=\tilde{C}$ in $NS(\tilde{Z})$. Indeed, $\frac{c}{t}$ is a rational function on $\tilde{Z}$ and $div(1-\frac{{\epsilon}\cdot c}{t})=Q'+L+E_1+E_2+E_3+E_4+E_5+E_6- \tilde{C}$.\\

Define $H=E_1+E_2+E_3+E_1'+E_2'+E_3'+Q=E_3+E_4+E_5+E_3'+E_4'+E_5'+Q=E_5+E_6+E_1+E_5'+E_6'+E_1'+Q=E_1+E_2+E_3+E_4+E_5+E_6+L\in NS(\tilde{Z})$. From the definition we see that $H^2=4$ and the complete linear system $\mid H \mid$ is base point free. Hence it gives a morphism ${\phi}{\mid}_{\mid H \mid}\colon \tilde{Z}\rightarrow {\mathbb P}^3$, which is an isomorphism, except that it contracts $E_1, E_3, E_5$ and $Q$ into $4$ nodes on ${\phi}{\mid}_{\mid H \mid}(\tilde{Z})$.\\  

Indeed, one can apply Theorem 5.2 of Saint-Donat \cite{SD} to check that $H$ is not hyperelliptic (see \cite{SD}). If it were, then there would be an irreducible curve $E\subset \tilde{Z}$ such that $E^2=0$ and $E\cdot H=2$. Since $H=r^{*}q^{*}h+Q$ and $Q\cdot E\geq 0$, $q(r(E))$ would be either a line (if $Q\cdot E=1$), or a conic (if $Q\cdot E=0$) in ${\mathbb P}^2$. Hence $E\rightarrow q(r(E))\cong {\mathbb P}^1$ would be a double cover ramified over $\leq E\cdot \tilde{C}=E\cdot H+E\cdot Q'=E\cdot H+E\cdot Q=deg(q(r(E)))+2E\cdot Q<4$ points. This is impossible.\\ 

Then the images of $E_2, E_4, E_6, E_1', E_3', E_5'$ are the lines in ${\mathbb P}^3$ connecting the nodes ${\phi}{\mid}_{\mid H \mid}(E_1)$, ${\phi}{\mid}_{\mid H \mid}(E_3)$, ${\phi}{\mid}_{\mid H \mid}(E_5)$, ${\phi}{\mid}_{\mid H \mid}(Q)$. So, if we denote $X={\phi}{\mid}_{\mid H \mid}(\tilde{Z})$, then $X$ is a $4$-nodal tetrahedral quartic in ${\mathbb P}^3$ and ${{\phi}{\mid}_{\mid H \mid}}{\mid}_X\colon \tilde{Z}\rightarrow X$ is its minimal resolution of singularities.\\

One also sees that the double cover $r\colon \tilde{Z}\rightarrow Z$ comes from the projection of the tetrahedral quzrtic $X$ from its node ${\phi}{\mid}_{\mid H \mid}(Q)$ and that $\tilde{C}$ is the normalization of the branch curve of this projection. {\it QED}\\

{\bf Remark.} Alternatively, one can describe curves $C\subset Z$ leading to tetrahedral quartics in ${\mathbb P}^3$ as follows (see \cite{Mukai}):

\begin{itemize}
\item[(*)] $C$ is a smooth genus $7$ curve, which is neither hyperelliptic, nor trigonal, nor bielliptic, has exactly $3$ $g^{1}_{4}$, has exactly $2$ $g^{2}_{6}$ (let us denote them by $\delta$ and ${\delta}'$) and contains $12$ distinct points $a_1,a_2,a_3,b_1,b_2,b_3,p_1,p_2,p_3,q_1,q_2,q_3$ such that $\delta = \mid 2a_1+2a_2+2a_3 \mid$, ${\delta}' = \mid 2b_1+2b_2+2b_3 \mid$, $dim \mid 2a_1+2a_2+2a_3-2p_i \mid \geq 1$ and $dim \mid 2b_1+2b_2+2b_3-2q_i \mid \geq 1$ for any $i$, and $a_1+a_2+a_3+b_1+b_2+b_3+p_1+p_2+p_3+q_1+q_2+q_3$ is a canonical divisor on $C$.
\end{itemize}

So, general tetrahedral quartics in ${\mathbb P}^3$ are exactly double covers of del Pezzo surfaces of degree $6$ ramified over a canonically embedded curve of genus $7$, which satisfies condition $(*)$ above.\\

{\bf Remark.} We have just seen above how to recover a general tetrahedral quartic $X\subset {\mathbb P}^3$ from the sextic $C(X)\subset {\mathbb P}^2$, which is the branch locus of the projection of $X$ from one of its nodes.\\

{\bf Theorem 4.4.} {\it Let us denote by $\nu \colon \tilde{C}(X)\rightarrow C(X)$ the normalization of the plane sextic $C(X)$, which is the branch curve of the projection of a general $4$-nodal tetrahedral quartic $X\subset {\mathbb P}^3$ from one of its nodes. Let ${\pi}_i\colon \tilde{X}_i\colon \rightarrow X_i$ ($i=1,2$) be the minimal resolution of singularities of two general $4$-nodal tetrahedral quartics in ${\mathbb P}^3$.\\

If $\tilde{C}(X_1)\cong \tilde{C}(X_2)$ as abstract genus $7$ curves, then $\tilde{X}_1\cong \tilde{X}_2$ as abstract $K3$ surfaces.}\\  

{\it Proof:} We will use the fact that $\tilde{C}(X)$ has exactly $2$ $g^{2}_{6}$ corresponding to two natural morphisms $\tilde{X}\rightarrow X\subset {\mathbb P}^3$ and $\tilde{X}\rightarrow X^{v}\subset {\mathbb P}^3$. Let us denote them by ${\delta}_i$ and ${{\delta}_i}'$ for the curve $\tilde{C}(X_i)$, $i=1,2$.\\ 

Let $\alpha \colon \tilde{C}(X_1){\rightarrow} \tilde{C}(X_2)$ be an isomorphism. Without loss of generality, we may assume that $(\alpha)^{*}({\delta}_2)={\delta}_1$. This implies that the images of $\tilde{C}(X_1)$ and $\tilde{C}(X_2)$ in ${\mathbb P}^2$ (under ${\delta}_1$ and ${\delta}_2$ respectively) are the same (upto a projective linear isomorphism), i.e. $C(X_1)=C(X_2)$ as plane sextics.\\

Now we can use the procedure above to recover $\tilde{X}_1$ and $\tilde{X}_2$ from $C(X_1)$ and $C(X_2)$ respectively, and we get that $\tilde{X}_1\cong \tilde{X}_2$, because they are double del Pezzo surfaces of degree $6$ ramified over the same canonically embedded curve of genus $7$. {\it QED}\\ 

\section{Quotients of $K3$ surfaces by an involution} 

We will show here how to view a general tetrahedral quartic $X$ as a (partial) resolution of singularities of a quotient of a $K3$ surface by a (Nikulin) involution with $8$ fixed points.\\

Let $H'=A-L_{14}-L_{23}$. Then $H'$ is nef, $H'\cdot E_i=H'\cdot R_j=0$ for all $i,j$, $H'\cdot L_{12}=H'\cdot L_{13}=H'\cdot L_{24}=H'\cdot L_{34}=2$, $H'\cdot L_{14}=H'\cdot L_{23}=4$, $(H')^2=8$.\\

It follows from criteria of \cite{SD} (see also \cite{Mori}, Theorem 5) that $H'$ is not hyperelliptic and the complete linear system $\mid H' \mid$ is base point free. Indeed, if $\mid H'\mid$ were not base point free, then by \cite{SD} $\mid H'\mid$ would have a fixed component $\Gamma$, which would be one of the curves $E_1, E_2, E_3, E_4, R_1, R_2, R_3, R_4$, $L_{12}, L_{13}, L_{24}, L_{34}$. By Theorem 5 in \cite{Mori} there would exist a smooth irreducible genus $1$ curve $E$ such that $H'=\Gamma + 5E$ and $\Gamma\cdot E=1$. This however implies that $H'\cdot \Gamma=3$, which is not the case for either of the hypothetical curves $\Gamma$ listed above. Hence $\mid H'\mid$ is base point free. In order to see that $H'$ is not hyperelliptic, let us apply Theorem 5.2 of \cite{SD}. If it were not, then one of the following two cases would occur. Either (case (i) in Theorem 5.2. of \cite{SD}) there would exist a smooth irreducible genus $1$ curve $E\subset \tilde{X}$ such that $E\cdot H'=2$, or (case (ii) in Theorem 5.2. of \cite{SD}) there would exist a smooth irreducible genus $2$ curve $B\subset \tilde{X}$ such that $H'=2B$. In the first case one arrives at a contradiction by considering an elliptic fibration ${\phi}_{\mid E\mid}\colon \tilde{X}\rightarrow {\mathbb P}^1$ corresponding to $E$ and noticing that at least two of the curves $L_{12}, L_{13}, L_{24}, L_{34}$ would have to lie in the same fiber of ${\phi}_{\mid E\mid}$, which would imply that, say, $E\geq L_{12}+L_{13}$, and so $2=H'\cdot E\geq H\cdot L_{12}+H\cdot L_{13}=4$, which is impossible. In the second case curve $B$ (its image under the minimal resolution of singularities ${\pi}\colon \tilde{X}\rightarrow X\subset {\mathbb P}^3$) would be a smooth nondegenerate genus $2$ and degree $4$ space curve, which is nonsense. So, $H'$ is not hyperelliptic and $\mid H' \mid$ is base point free.\\

Hence the complete linear system $\mid H' \mid$ determines a morphism ${\phi}{\mid}_{\mid H' \mid}\colon \tilde{X}\rightarrow {\mathbb P}^5$, which is the minimal resolution of singularities of its image ${\phi}{\mid}_{\mid H' \mid}(\tilde{X})$, which we will denote by $X'={\phi}{\mid}_{\mid H' \mid}(\tilde{X})$. Then $X'$ is a normal surface of degree $8$ in ${\mathbb P}^5$ with exactly $8$ nodes (which are the images of $E_1, E_2, E_3, E_4, R_1, R_2, R_3, R_4$). The images of $L_{13}, L_{14}, L_{23}, L_{24}$ form a cycle of $4$ conics in ${\mathbb P}^5$ (which is a hyperplane section of $X'$). The images of $L_{12}$ and $L_{34}$ are smooth rational normal (i.e. span a hyperplane) quartic curves in ${\mathbb P}^5$. 

\begin{center}
  \includegraphics{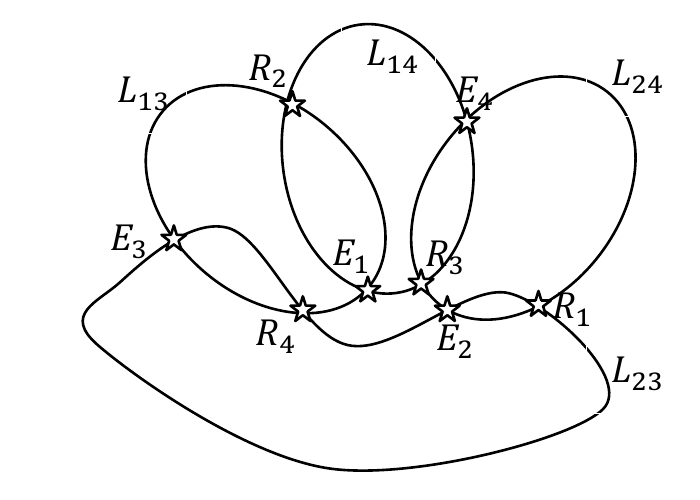}
\end{center}

Let ${\Pi}_{13}, {\Pi}_{14}, {\Pi}_{23}, {\Pi}_{24}$ be the $2$-planes in ${\mathbb P}^5$, spanned by the conics $L_{13}, L_{14}, L_{23}, L_{24}$ respectively. Then ${\Pi}_{13}$ and ${\Pi}_{14}$ intersect along a line $\overline{E_1R_2}$, ${\Pi}_{14}$ and ${\Pi}_{24}$ intersect along a line $\overline{E_4R_3}$, ${\Pi}_{23}$ and ${\Pi}_{24}$ intersect along a line $\overline{E_2R_1}$, and ${\Pi}_{23}$ and ${\Pi}_{13}$ intersect along a line $\overline{E_3R_4}$. One also notices that all these $4$ $2$-planes lie in a hyperplane in ${\mathbb P}^5$ and intersect at exactly one point $P$ (i.e. form a bouquet of $4$ $2$-planes), which is the point of intersection of the $4$ lines listed above.\\

\begin{center}
  \includegraphics{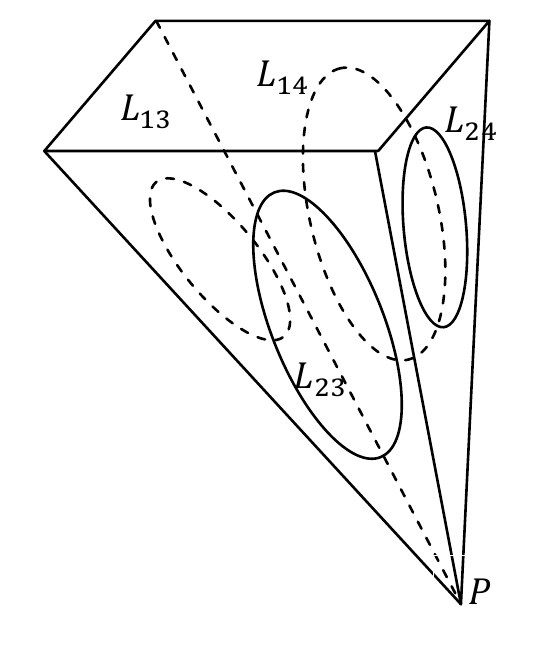}
\end{center}

Morphism $\tilde{X}\rightarrow X'$ factors through $X$ and represents it as a partial resolution of singularities of $X'$.\\

{\bf Lemma 5.1.} {\it If $X$ is a $4$-nodal tetrahedral quartic, then $E_1, E_2, E_3, E_4, R_1, R_2, R_3, R_4$ form an even eight on its minimal resolution of singularities $\tilde{X}$.}\\

{\it Proof:} $E_1+E_2+E_3+E_4+R_1+R_2+R_3+R_4=2\cdot (L_{23}+L_{24}+L_{34}-L_{12}-L_{13}-L_{14}-E_1+E_2+E_3+E_4+2R_1)$ in $M\subset NS(\tilde{X})$. {\it QED}\\

Let $\theta \colon Y\rightarrow \tilde{X}$ be the double cover of $\tilde{X}$ ramified over $8$ disjoint curves $E_1, E_2$, $E_3, E_4$, $R_1, R_2$, $R_3, R_4$. It is immediate that $Y$ is a blowup of a $K3$ surface and $\tilde{X}$ is its quotient by the covering involution $i\colon Y {\rightarrow} Y$. Let $E_1', E_2', E_3', E_4', R_1', R_2', R_3', R_4'$ be the reduced preimages of the branch curves. Notice that each of them is a $(-1)$-curve, and so there is a contraction morphism ${\pi}'\colon Y\rightarrow Y'$, which contracts $8$ curves $E_1', E_2', E_3', E_4'$, $R_1', R_2', R_3', R_4'$ to $8$ fixed points of the involution $i'\colon Y'\rightarrow Y'$ induced on $K3$ surface $Y'$ by $i$.\\

One has the following commutative diagram:\\

\begin{center}
  \includegraphics{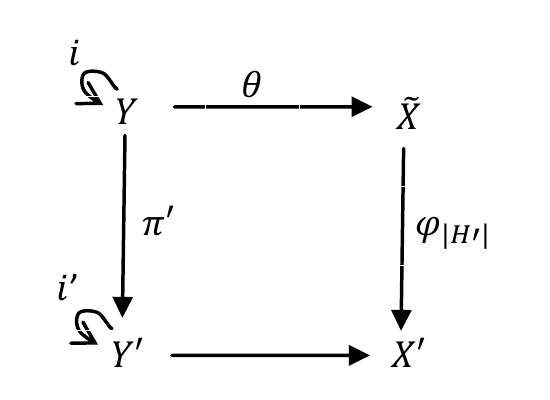}
\end{center}

which represents $X'$ as a quotient of a $K3$ surface $Y'$ by an involution $i'$ with exactly $8$ fixed points. This involution is symplectic, because $X'$ is a $K3$ surface by construction. This proves:\\

{\bf Theorem 5.2.} {\it A general tetrahedral quartic in ${\mathbb P}^3$ is a partial resolution of singularities of a quotient of a $K3$ surface by a Nikulin involution with $8$ fixed points.}\\

{\bf Question.} How to characterize $K3$ surfaces $Y'$ (or maybe pairs $(Y',i')$), which lead to tetrahedral quartics in ${\mathbb P}^3$ by taking their quotient by a symplectic involution with $8$ fixed points?\\

\section{Acknowledgement}

We thank Professor Dolgachev for finding and sending to us Bauer's article \cite{Bauer}.\\

\bibliographystyle{ams-plain}

\bibliography{TetrahedralQuartics}

\end{document}